\documentclass[journal]{IEEEtran}

\usepackage{pgfplots}
\usepackage{pgfplotstable} 
\usepackage{filecontents}

\usepackage{soul,framed} 
\usepackage{xcolor}
\colorlet{shadecolor}{yellow}
\DeclareGraphicsExtensions{.pdf,.jpeg,.png}
\usepackage{verbatim}
\usepackage[cmex10]{amsmath}
\usepackage{array}
\usepackage{mdwmath}
\usepackage{eqparbox}
\usepackage{url}
\usepackage{amsfonts,amssymb}
\usepackage{amsthm}
\usepackage{fontawesome}
\usepackage{tikz,pgfplots}
\usepackage{multirow}
\usepackage{footnote}
\makesavenoteenv{tabular}
\makesavenoteenv{table}
\usepackage{booktabs}

\theoremstyle{definition}
\newtheorem{definition}{Definition}
\newtheorem{rmk}{Remark}
\theoremstyle{plain}

\def\bitcoinA{%
	\leavevmode
	\vtop{\offinterlineskip 
		\setbox0=\hbox{B}%
		\setbox2=\hbox to\wd0{\hfil\hskip-.03em
			\vrule height .3ex width .15ex\hskip .08em
			\vrule height .3ex width .15ex\hfil}
		\vbox{\copy2\box0}\box2}}

\newcommand{\expected}{\mathop{\mathbb{E}}}
\def\tablename{Table}
\pgfplotsset{%
	compat=newest, 
	tick label style={font=\footnotesize},
	label style={font=\footnotesize},
	legend style={font=\footnotesize},
	axis x line = center,
	axis y line = center,
	every axis/.style={pin distance=1ex},
	trim axis left
}

\usepgfplotslibrary{dateplot}
\usepgfplotslibrary{fillbetween}
\usepgfplotslibrary{colorbrewer}
\usetikzlibrary{shapes.multipart}
\usetikzlibrary{automata,positioning}

\newcommand{\subparagraph}{}

\begin{document}
	
	\title{Integrated Planning of a Solar/Storage Collective
		\thanks{This work was funded by the Make Our Planet Great Again Program of the French Government and the interdisciplinary Eco-SESA Program of the Univ. Grenoble Alpes.}
	}
	
	\author{Jesus E. Contreras-Oca\~na, Arshpreet Singh, Yvon B\'esanger, Fr\'ed\'eric Wurtz}

	
	\maketitle

	\begin{abstract}
		French regulation allows consumers in low-voltage networks to form collectives to produce, share, and consume local energy under the \textit{collective self-consumption framework}. A natural consequence of collectively-owned generation projects is the need to allocate production among consumers. In long-term plans, production allocation determines each of the consumers' benefits of joining the collective. In the short-term, energy should be dynamically allocated to reflect operation. This paper presents a framework that integrates long and short-term planning of a collective that shares a solar plus energy storage system.  In the long-term planning stage, we maximize the collective's welfare and equitably allocate expected energy to each consumer. For operation, we propose a model predictive control algorithm that minimizes short-term costs and allocates energy to each consumer on a 30-minute basis (as required by French regulation). We adjust the energy allotment ex-post operation to reflect the materialization of uncertainty. We present a case study where we showcase the framework for a 15 consumer collective.      
	\end{abstract}

	\begin{IEEEkeywords}
		Planning, control, energy community, techno-economic modeling. 
	\end{IEEEkeywords}

	\IEEEpeerreviewmaketitle
	\section*{Nomenclature}
	\subsection*{Indices and sets}
	\begin{IEEEdescription}[\IEEEusemathlabelsep\IEEEsetlabelwidth{$xxxx$}]
		\item[$i$] Index of consumers
		\item[$\mathcal B$] Set of feasible ES charge/discharge actions
		\item [$\mathcal E$] Set of feasible annual energy allocation per consumer
		\item[$\omega$] Index of uncertainty scenarios
		\item[$\Omega$] Set of uncertainty scenarios
	\end{IEEEdescription}
	
	\subsection*{Variables and parameters}
	\begin{IEEEdescription}[\IEEEusemathlabelsep\IEEEsetlabelwidth{$xxxx$}]
		\item[$A$] Number of years in planning horizon
		\item[$\overline B_\mathrm{net}$] Expected net benefit (EUR)
		\item[$c/d$] ES Charged/discharged per time period (kWh)
		\item[$C_i$] Consumer $i$'s cost (EUR)\footnote{$l$, $g$, and $C$ (without subscripts) represent aggregate load, local  generation, and cost, respectively. }
		\item[$C_\mathrm{grid}$] Grid connection cost (EUR)
		\item[$C^\mathrm{es}_\mathrm{inv}$] ES inverter cost (EUR)
		\item[$C^\mathrm{pv}_\mathrm{inv}$] PV inverter cost (EUR)
		\item[$\mathrm{CapEx}$] Capital expenditures (EUR)
		\item[$e$] Energy allocation per consumer (kWh) \footnote{See Section~\ref{subsec:Control} for the meaning of each of the superscripts and accents of $e$. }
		\item[$E^\mathrm{es}_\mathrm{cap}$] Energy storage capacity (kWh)
		\item[$G$] Key of repartition matrix (kWh)
		\item[$g_i$] Local generation assigned to consumer $i$ (kWh)
		\item[$g^\mathrm{g}$] Energy bought from the grid per time period (kWh)
		\item[$g^\mathrm{pv}$] PV generation per time period (kWh)
		\item[$g^\mathrm{s}$] Surplus energy sold to the grid per time period (kWh)
		\item[$L$] Consumer load per time period matrix (kWh)
		\item[$l_i$] Load of consumer $i$ (kWh)
		\item[$\overline m$] Expected mismatch between planned and delivered energy (kWh)
		\item[$N$] Number of consumers
		\item[$\mathrm{OpEx}$] Operational expenditures (EUR)
		\item[$\mathrm{obj}^*_\mathrm{LT}$] Objective of the long-term planning problem (EUR)
		\item[$p$] Price of local PV+S energy (EUR/kWh)
		\item[$P^\mathrm{es}_\mathrm{cap}$] Energy storage power capacity (kW)
		\item[$P^\mathrm{es}_\mathrm{inv}$] Energy storage inverter capacity (kW)
		\item[$P^\mathrm{pv}_\mathrm{cap}$] PV installed capacity (kW)
		\item[$P^\mathrm{pv}_\mathrm{inv}$] PV inverter capacity (kW)
		\item[$r$] Discount rate
		\item[$S$] Subsidies (EUR)
		\item[$T$] Number of time periods for simulation in a year\footnote{The subscripts $c$ and $p$ denote control and prediction horizon, respectively.} 
		\item[$\alpha$] Vector of per-unit PV production
		\item[$\beta^\mathrm{es}/\beta^\mathrm{pv}$] Per-unit cost of ES/PV capacity (EUR/kWh)/(EUR/kW)
		\item[$\beta^\mathrm{es \;u}$] ES utilization cost coefficient (EUR/kWh)
		\item[$\beta^\mathrm{mnt}$] Yearly PV maintenance cost per unit of intalled capacity (EUR/kW)
		\item[$\Delta$] Length of time period (hours)
		\item[$\hat \epsilon$] Deviation between actual and expected PV+S generation in the control horizon (kWh) 
		\item[$\gamma$] Share of the net benefit kept by investor 
		\item[$\kappa$] Energy-to-power ratio of the energy storage system
		\item[$\lambda$] Price of energy sold to the grid (EUR/kWh)
		\item[$\Pi$] Investor's profit (EUR)
		\item[$\tau$] Grid export tax (EUR/kWh)
		\item[$\theta$] Weight of energy mismatch (EUR/kWh)
	\end{IEEEdescription}
	
	\section{Introduction}
	
	As part of an almost universal embrace of renewable energies by people and governments around the globe, many jurisdictions today recognize and encourage the \textit{self-consumption} of local renewable energy generation. For example, the European Commission put out a set of best practices for self-consumption regulation~\cite{eu2015best} and several European countries including Germany, France, Spain, and Switzerland have in place laws and frameworks related to self-consumption~\cite{Frieden_2019_collective}.
	
	In France, the practice of self-consumption (\textit{autoconsommation}) was recognized by the energy code in 2017~\cite{code_energie}. The code defines the rate of self-consumption as the share of the local production at each instant\footnote{In practice and as dictated by France's Energy Regulatory Commission, the rate of self-consumption is calculated every 30 minutes~\cite{deliberation_CRE}.} that is consumed locally.  Additionally, the code lays out a framework for \textit{collective} self-consumption. In the collective variant, a group of geographically proximal consumers connected to the low-voltage grid can share locally produced energy~\cite{code_energie}. In the French case, collectives organize around a legal entity (\textit{personne morale}) who is in charge of communicating the repartition of energy among consumers on a 30-minute basis. 
	
	From the consumer's point of view, we identify two major advantages of collective over individual self-consumption. First, it facilitates investment by increasing the viability of larger systems and enable economies of scale~\cite{abada2017viability}, opens up otherwise unavailable physical resources (e.g., better-located rooftops), and allows for the participation of institutional investors that can provide attractive financing and absorb risk~\cite{Dunlop_2016_EU, goedkoop2016partnership}. The second advantage is that by sharing energy, collectives can reach higher rates of self-consumption than equivalent individual projects. This is especially true if the loads are heterogeneous. Thus, if the electricity tariff structure rewards self-consumption, collectives reap higher benefits.        
	
	However, organizing a collective also represents challenges. For one, the interconnection process and operation one large system be complicated~\cite{shared_schroeder}. Or, if one aggregates several smaller systems, then the challenge becomes one of coordination and control. In any case, coordination among consumers and the grid represents a challenge~\cite{Council_2019_regulatory}. Beyond technical challenges, several studies identify the difficulty of reaching an acceptable deal to all parties - investor and consumers - as significant challenges~\cite{abada2017viability, goedkoop2016partnership, Council_2019_regulatory}.  In particular, Goedkoop and Devine-Wright in~\cite{goedkoop2016partnership} stress the need for mechanisms that create trust among the participants and lead to stable collective arrangements.  
	
	\subsection{Summary and contribution of this work}

	In this paper, we present \emph{a cooperative planning framework that integrates long-term planning and short-term operation of an energy collective} of consumers sharing a photovoltaic plus storage (PV+S) system and an investor that provides capital. By ``cooperative'', we mean a planning framework that puts the welfare of the collective ahead of individuals. We designed the framework with cooperation rather than competition as a guiding principle because we agree with Goedkoop and Devine-Wright: an overly-competitive environment may lead to tense relations, lack of trust, and ultimately jeopardize the financial stability of the collective~\cite{goedkoop2016partnership}. 
	
	Furthermore, we adopt simplicity\footnote{Developing a technical and quantitative meaning of ``simple'' is not trivial and outside the scope of our work. In our work, mechanism 1 is simpler than mechanism 2 if we believe that it is easier to explain in layman's terms and its input-output relationship is easier to understand.} as a guiding design principle for the mechanisms that frame the relationship between participants. As outlined in~\cite{eu2015best} simple procedures are important if we want to encourage healthy energy collectives. Concretely, our framework is designed with transparent and simple-to-understand prices and operational rules.    
	
	In our framework, the long-term plan answers two important questions: \emph{what PV+S system size is optimal? and what benefits can consumers and the investor expect from joining the collective?} The answer to the first question is important because we want to use resources as efficiently as possible. Knowing how much each consumer can expect to save and how much the investor can expect to earn is important because it allows them to decide whether or not to join the collective. 
	
	The main determinant of the investor's profit is the price $p$ at which PV+S energy is sold to consumers. For the investor, $p$ should be high enough to recover its capital and operating costs (CapEx and OpEx). However, for the consumers to benefit from joining the collective, $p$ should be (roughly speaking) lower than grid prices. In this work, we \emph{present a method to determine the range of values that $p$ can take such that all participants benefit.} The existence of such a range is necessary for the collective to be \textit{financial sustainable}, i.e., a collective in which all participants benefit financially.
	
	\begin{figure}
\includegraphics[width=0.4\linewidth]{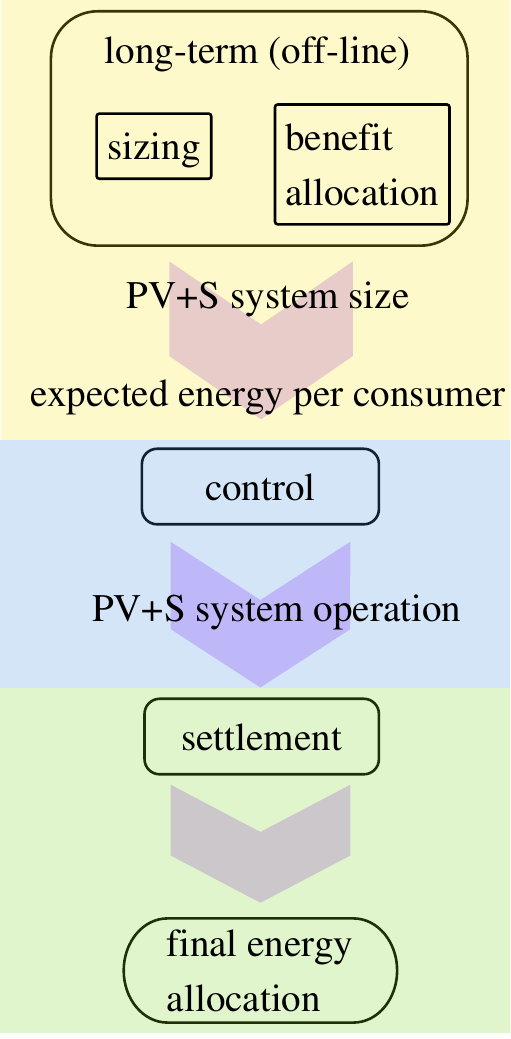}
\centering
\caption{Diagram representing the framework proposed in this paper.}
\label{fig: outline}
\end{figure}

	Part of the long-term plan is to \textit{calculate the repartition of local energy production}. French regulation requires the legal entity to produce a ``key of repartition'' to determine the share of local production received by each consumer on a 30 minutes basis. Thus, we propose a \textit{method to determine a key of repartition that equitably allocates energy among consumers throughout a year}. This key determines the expected amount of PV+S energy that each consumer can expect to receive in a year. The expected energy allotment for each consumer and the price $p$ are the two ingredients needed for each consumer to estimate the benefits of joining the collective. The yellow section of Fig.~\ref{fig: outline} illustrates the long-term planning portion of our framework.
	
	The two main goals of the real-time operation portion of our framework are to efficiently operate the PV+S system and fulfill the energy expectations of each customer. Real-time operation is composed of a control stage (the middle blue section in Fig.~\ref{fig: outline}) and a settlement stage (the green section of Fig.~\ref{fig: outline}). In the control stage, we \textit{propose a model predictive control (MPC)- based algorithm to minimize operating costs and determine a key of repartition that tracks the expected energy allotment from the long-term plan.} In the settlement stage, we \textit{propose an optimization-based method to adjust the key of repartition from the control strategy} such that it complies with the actual operation (which differs from the control strategy due to materialization of uncertainty) of the PV+S system. 
	
	We would like to emphasize that energy allocation is done for financial and billing purposes. That is, for each consumer, the algorithm decides how much of their energy consumption should be paid to the grid and how much to the local PV+S. We do not control the actual power flows as these are determined by Kirchhoff's Laws.
	
	We would like to emphasize that while this work was inspired by the French collective self-consumption case, the proposed framework is generic and widely applicable. Within European regulation, the proposed framework could be used under the definitions for ``Renewable Energy Communities'' from~\cite{EU_RE_directive} and ``Citizen Energy Communities'' from~\cite{EU_IME_directive}. To this day, Portugal, Spain, Wallonia (Belgium), and Greece have adopted regulations on energy communities and the rest of the European Member States should follow within the next few years. Beyond European regulation, our framework could be used for the planning and operation of a generic microgrid with community or investor-owned assets.  
	
	Finally, we provide a \textit{case study of a PV+S energy collective} in the South of France. In it, we showcase each of the elements of the proposed integrated planning framework using realistic data.
	
	To recapitulate the contributions of this work, we restate them as follows. 
	\begin{itemize}
		\item A cooperative planning framework that integrates long-term planning and short-term operation of an energy collective.
		\item A long-term plan determines that determines the optimal PV+S system and the conditions for the financial sustainability of the collective.
		\item A method to equitably allocate yearly local PV+S energy to the consumers. 
		\item An MPC-based algorithm to minimize operating costs and determine the key of repartition that tracks the expected energy allotment of the long-term plan.
		\item A method to adjust the key of repartition ex-post operation to reflect the realization of uncertainties and operation of the system.
		\item A case study of our framework for a 15 consumer energy collective.
	\end{itemize}

	\subsection{Literature review}
	
	Energy communities are not an entirely new concept (long-existing cooperatives in the U.S. and microgrids are related concepts) but following the unveiling of the Clean Energy Package for All Europeans set of legislation, have gained renewed interest~\cite{Frieden_2018_Overview, alaton2020energy}. In~\cite{Moret_2019_Energy}, Moret and Pinson introduce the concept of an ``Energy Collective" composed of prosumers that interact with each other and with the grid through a third party. Their work introduces a structure similar to the one presented in this paper and analyses the operation and fairness of the collective. However, it does not explicitly deal with the problem of energy repartition by community-owned assets and ignores the problem of long-term planning. The work in~\cite{Vanadzina_2019_Business} presents a review of several business models under which a collective can be organized. Some of those business models could in principle be applied within our framework (e.g., customer and mixed ownership, energy service company).  
	
	One of the central issues introduced by the collective self-consumption framework is the problem of local energy allocation among consumers. Market-based mechanisms are some of the most common ways of allocating shared resources, e.g., the works in~\cite{Fleischhacker_2019_Sharing, Nunna_2017_Multiagent, Park_2017_Event, Shamshi_2016_Shamsi, Liu_2019_ANovel, Feng_2019_Caolitional} for the case of microgrids. Alternative non-market methods of resource allocation also exist. For example, in~\cite{Etinski_2013_Fair}, Erinski and Sch\"ulke propose a method to fairly allocate energy from a central wind or solar generator. In~\cite{Chen_2019_PV}, the authors propose an energy credit-based system to share solar energy surplus within a community. The authors of~\cite{abada2017viability} and~\cite{Feng_2019_Caolitional} present cooperative game-theoretic analyses of energy collectives. In~\cite{Feng_2019_Caolitional}, the authors present a nucleolus-based solution that leads to a stable and fair payoff distribution
	scheme for all players. Abada et al.'s work in~\cite{abada2017viability} touches on many of the questions that we deal with in this paper, e.g., optimal investment and financial stability of an energy collective. While~\cite{abada2017viability} presents a deeper and more theoretical analysis of the question of financial stability, they ignore some of the elements of our work such as the energy storage (ES) system, the key of repartition, real-time operation, and ex-post operation settlements.

	PV+S system sizing is a key topic of this paper and also a heavily studied subject. Works such as as~\cite{Khare_2014_Optimal, Kolhe_2009_Techno, Contreras-Ocana_2019_Non-Wire,Shahidehpour_2019_Optimal} focus on the techno-economics of sizing while works such as as~\cite{Xu_2013_AnImporved, Akram_2018_AnImproved, Saez-de-Ibarra_2019_Co-Optimization} direct their focus to more technical aspects. In~\cite{Khare_2014_Optimal}, Khare and Ragnekar develop a methodology that minimizes the annual cost of a PV system. The work in~\cite{Kolhe_2009_Techno} goes a step further and considers a PV+S system. The authors of~\cite{Contreras-Ocana_2019_Non-Wire} present a sizing problem of distributed energy resources that incorporates the value of delaying capacity expansion. The work in~\cite{Shahidehpour_2019_Optimal} presents a planning problem of an isolated microgrid. The algorithm in~\cite{Shahidehpour_2019_Optimal} is especially interesting because it exploits the complementary of solar energy and biogas to reduce the need for ES. The authors of~\cite{Xu_2013_AnImporved} propose a sizing method for solar, wind, and battery microgrids that considers reliability metrics, the specific characteristics of wind and solar production, among others. The work in~\cite{Akram_2018_AnImproved} proposes a microgrid sizing algorithm that, in addition to wind, solar, and batteries, also considers dispatchable diesel generation. Similar to this work,~\cite{Contreras-Ocana_2019_Non-Wire} and~\cite{Saez-de-Ibarra_2019_Co-Optimization}, embed the short-term operation of the system in their sizing problems.  However, the question of how to allocate local production is outside the scope of all the aforementioned works.

	\subsection{Organization of the rest of the paper}
	
	Section~\ref{sec:participants} describes the consumer, investor, and legal entity models. Section~\ref{sec:long-term planning} presents the PV+S system sizing problem and Section~\ref{sec:benefit} presents the problem of benefit allocation. Section~\ref{sec:RTC} presents the system operation algorithms (control and ex-post operation settlement). Section~\ref{sec:casestudy} presents a case study of a PV+S collective in the south of France and Section~\ref{sec:Conclusions} concludes this paper and provides suggestions for future work. 
	
	\section{Participants in the Collective}
	\label{sec:participants}
	\subsection{Consumers}
	
	We model consumers as inflexible loads partially fulfilled by local energy and partially by grid energy. Let $l_i$ be a vector $\in \mathbb R_+^T$ that represents energy consumption of consumer $i$ during $T$ periods of length $\Delta$. Let $g_i$ be a vector $\in \mathbb R_+^T$ that represents energy from the PV+S system assigned to $i$ during each period. Then, $i$'s costs are
	\begin{equation}
	C_i( l_i - g_i) + p\cdot  1^\intercal g_i. \label{eq:resident_cost}
	\end{equation}
	where the first term is grid costs and the second local energy costs. The linear function $C_i$ maps $i$'s energy deficit, $l_i - g_i$, to grid-related cost which includes energy, fees, and taxes. The price of local energy is represented by $p$.   
	
	We express the total cost, $C$, of $N$ consumers as a function of the aggregate consumption and generation, $l - g = \sum_{i=1}^N l_i - g_i$:
	\begin{equation}
	C( l - g) +  p \cdot 1^\intercal g =  \sum_{i=1}^N \left\{C_i( l_i - g_i) +  p\cdot 1^\intercal g_i \right \}. \label{eq:all_resident_cost}
	\end{equation}
	
	\subsection{Investor}
	
	The investor provides capital for the PV+S system, funds maintenance, operational expenses, taxes, and absorbs the bulk of the risk\footnote{The main sources of risk are solar generation uncertainty, future prices, and regulatory uncertainty. Consumers are also exposed to risk since solar supply is uncertain.}. In return, the investor profits by selling PV+S to the consumers. 
	
	We define the present value of the investor's profit as 
	We define the present value of the investor's profit as 
	\begin{equation}
	\Pi = \underbrace{\sum_{a=1}^A \frac{p \cdot 1^\intercal g + \lambda^\intercal g^\mathrm{s}  + S }{(1+r)^a} }_\text{revenues} - \underbrace{\mathrm{CapEx} - \sum_{a=1}^A \frac{\mathrm{OpEx}}{(1+r)^a}}_\text{costs}. \label{eq:investor_profit}
	\end{equation}
	The investor's revenue stream is composed of three elements: sales to the consumers, $p\cdot  1^\intercal g$, energy sales to the grid\footnote{Energy sales to the grid can be limited by technical or regulatory limits.} at a price $\lambda$, $ \lambda^\intercal g^\mathrm{s}$, and yearly subsidies, $S$. To calculate the present value of the total revenue, we add the discounted revenue stream over the planning horizon $\{1,2,\hdots,A\}$. Here, $r$ is the discount rate and relates future cash flows to present value~\cite{levy_1994_capital}.   
	
	The investor's costs are composed of CapEx and OpEx. The former are incurred before and during system installation and include hardware, equipment, installation, and grid connection costs. The latter include system utilization costs and taxes and, similar to the revenue stream, we add the discounted OpEx stream to calculate the present value of the OpEx. 
	
	\subsection{Legal entity}

	The third participant is a neutral legal entity who is in charge of operating the PV+S system, allotting local energy to consumers, billing, and representing the collective. It behaves neutrally and is guided by contracts and algorithms. Concretely, the legal entity determines and implements the PV+S system real-time operation strategy, calculates the ex-post operation settlements, bills, and communicates with the grid.

	The EU Directive 2019/944 on common rules for the internal market for electricity and the EU Directive 2018/2001 on the promotion of the use of energy from renewable sources define ``Citizen Energy Communities'' and  ``Renewable Energy Communities'' as legal entities (e.g., civil associations, cooperatives, corporations, etc.). For example, France's legislation calls for self-consumption activities to be organized within a \textit{personne morale} (i.e., a legal entity or legal person)~\cite{code_energie}.  In Greece, natural persons and small and medium enterprises are allowed to create for-profit or non-for-profit cooperatives to legally house an energy community~\cite{Hannoset_2019_Energy}. These are just a few examples of legislation that mandates collectives to operate within a legal entity. However, in our work, we adopt a wider definition that may include, for example, aggregators or contracts between individuals. Despite its name, more than a legalistic definition, the key is that the legal entity undertakes the responsibilities listed in the previous paragraph.

	\section{Long-term planning (sizing problem)}
	\label{sec:long-term planning}
	
	In this work, we take a cooperative approach to long-term planning. That is, we look for a Pareto-optimal PV+S system that maximizes the total welfare of the investor and the consumers, i.e., the investor's profit minus the consumers' costs. 
	
	We formulate the sizing problem as a two-stage stochastic mixed-integer linear program (MILP) whose main source of uncertainty is solar irradiation. Solar irradiation uncertainty is modeled using scenarios of per-unit of installed PV capacity PV generation. Each scenario is a possible  PV generation profile and is represented by a $T$-long time series of numbers between 0 (e.g., at night) and 1 (i.e., generating at capacity).  
	
	The first stage (here-and-now) variables, i.e., decisions taken before the uncertainty materializes, are sizing-related decisions such as the capacities of the PV and ES systems. On the other hand, the second stage (wait-and-see) variables are decisions that can be modified later when refined forecasts are available. Specifically, second stage variables are operation-level decisions such as ES operation and PV production.   
	
	\subsection{Objective}
	
	We cast the problem as a maximization problem whose objective, 
	\begin{equation}
	\mathrm{obj}^*_\mathrm{LT}= \max 
	\left\{ \expected \left[\Pi_\omega \right]  - \sum_{a=1}^A { \frac{ \expected \left[ C( l - g_{\omega}) +p \cdot \!1^\intercal g_{\omega}\right] }{(1+r)^a} }\right\},  \label{eq:global_problem}
	\end{equation}
	is to maximize the expected welfare of the investor and the consumers. In this work, the subscript $\omega$ denotes scenarios of solar generation and the function $\expected$ calculates the expected value of its argument. It is trivial to show that the term $p\cdot 1^\intercal g_\omega$ cancels out in~\eqref{eq:global_problem} which renders the objective as linear.
	
	Our long-term plan simulates an entire year (rather than typical days as it is often the case in planning formulations). Since the OpEx (contained in $\Pi_\omega$) and grid costs depend on the short-term (e.g., hourly) operation of the PV+S system, we co-optimize investment and operation decisions. Including operation-scale decisions and scenarios of uncertainty makes the resulting MILP relatively large. Thus, to ease the computational burden in our case study, we model one typical year and assume that the rest of the planning horizon is comparable. The rest of this section describes the sizing and operational constraints of the problem. 
	
	\subsection{Sizing constraints}
	\label{subsec:Sizing constraints}
	The sizing-related constraints of the PV system are
	\begin{subequations}
		\begin{align}
		& (P^\mathrm{pv}_\mathrm{inv}, C^\mathrm{pv}_\mathrm{inv}) \in \{ (a_1^\mathrm{pv} , b_1^\mathrm{pv} ),(a_2^\mathrm{pv} , b_2^\mathrm{pv} ) \hdots \} \label{eq:PV_inverter_choice}\\  
		& P^\mathrm{pv}_\mathrm{cap} \le P^\mathrm{pv}_\mathrm{inv}. \label{eq:PV_CAP_LIM}
		\end{align}
	\end{subequations}
	Eq.~\eqref{eq:PV_inverter_choice} states that one must choose a PV inverter with maximum power $P^\mathrm{pv}_\mathrm{inv}$ and cost $C^\mathrm{pv}_\mathrm{inv}$ among a discrete set of inverters with characteristics $(a_j^\mathrm{pv} , b_j^\mathrm{pv})$ where $a_j^\mathrm{pv}$ is the maximum power of inverter $j$ and $b_j^\mathrm{pv}$ is its cost. Eq.~\eqref{eq:PV_CAP_LIM} limits PV capacity to the inverter's maximum power. We implement Eq.~\eqref{eq:PV_inverter_choice} in an MILP setting using discrete variables and linear constrains.
	
	The constraints of ES sizing decisions are
	\begin{subequations}
		\begin{align}
		& (P^\mathrm{es}_\mathrm{inv}, C^\mathrm{es}_\mathrm{inv}) \in \{ (a_1^\mathrm{es} , b_1^\mathrm{es} ),(a_2^\mathrm{es} , b_2^\mathrm{es} ) \hdots \} \label{eq:ES_inverter_choice}\\  
		& P^\mathrm{es}_\mathrm{cap} \le P^\mathrm{es}_\mathrm{inv} \label{eq:ES_CAP_LIM} \\ 
		& E^\mathrm{es}_\mathrm{cap} = \kappa \cdot P^\mathrm{es}_\mathrm{cap}. \label{eq:E_P_ES_proportion}
		\end{align}
	\end{subequations}
	Similar to the PV case, Eq.~\eqref{eq:ES_inverter_choice} defines the capacities and costs of the possible ES system inverters and Eq.~\eqref{eq:ES_CAP_LIM} limits the ES maximum charge/discharge rate to the inverter capacity. Eq.~\eqref{eq:E_P_ES_proportion} relates the ES system energy capacity to its maximum power via a fixed energy-to-power ratio, $\kappa$. 
	
	\subsection{Short-term operation constraints }
	
	Let the vector $g^\mathrm{pv}_\omega \in \mathbb R_+^T$ denote solar production in uncertainty scenario $\omega$ during each period. For each scenario, $g^\mathrm{pv}_\omega$ is the product of installed PV capacity, the length of each period $\Delta$, and a parameter $\alpha_\omega\in \mathbb R_+^T$ that is related to solar irradiation:
	\begin{equation}
	g^\mathrm{pv}_\omega =\Delta \cdot \alpha_\omega \cdot P^\mathrm{pv}_\mathrm{cap}  \; \forall \; \omega \in \Omega. \label{eq:def_g_pv}
	\end{equation}
	Here, $\Omega$ is the set of solar generation scenarios. Each of the elements of $\alpha_\omega$ takes values in $[0,1]$. When an element of $\alpha_\omega$ is $0$, there is no solar production and when it is $1$, the PV system produces at full capacity. 
	
	The ES system has two main operational variables, energy discharged and energy charged at each time period, denoted by $d_\omega\in \mathbb R_+^T$ and by $c_\omega\in \mathbb R_+^T$, respectively. In our work, we consider charge, discharge, and state-of-charge limits and charge and discharge efficiency. For the sake of brevity, we omit the details of the widely-used ES model (e.g., as in~\cite{Contreras-Ocana_2019_Non-Wire, Xu_2013_AnImporved, Akram_2018_AnImproved, Saez-de-Ibarra_2019_Co-Optimization, Contreras_2019_Participation}) and compactly denote these constrains as
	\begin{equation}
	(d_\omega, c_\omega) \in \mathcal B(P^\mathrm{es}_\mathrm{cap},E^\mathrm{es}_\mathrm{cap}) \; \forall \; \omega \in \Omega \label{eq:def_ES_const}
	\end{equation}
	where $\mathcal B(P^\mathrm{es}_\mathrm{cap},E^\mathrm{es}_\mathrm{cap})$ is the feasible charge and discharge space as a function of the battery size.

	We define imports from the grid, $g^\mathrm{g}_\omega\in \mathbb R_+^T$, and surplus of the collective, $g^\mathrm{s}_\omega\in \mathbb R_+^T$ as follows
	\begin{subequations}
		\begin{align}
		& g^\mathrm{g}_\omega =  \left[l_\omega  + c_\omega -g^\mathrm{pv}_\omega  - d_\omega   \right ]^+ \; \forall \; \omega \in \Omega \label{eq:grid_energy_def} \\
		& g^\mathrm{s}_\omega =  \left[g^\mathrm{pv}_\omega  + d_\omega - c_\omega  -l_\omega \right]^+ \; \forall \; \omega \in \Omega. \label{eq:surplus_energy_def}
		\end{align} \label{eq:energy_direction_def}
	\end{subequations} 
	We define imports from the grid as the difference of consumption (the sum of load and battery charge, $l_\omega + c_\omega$) and production (the sum of PV generation and battery discharge, $g^\mathrm{pv}_\omega  + d_\omega $) when this amount is positive. The operator $[ \;]^+$ is an element-wise function that gives the positive component of its argument. Similarly, Eq.~\eqref{eq:surplus_energy_def} defines surplus as the difference between production and consumption when this number is positive. 
	
	Another relevant quantity is PV+S generation assigned to the consumers, $g_\omega$, during each period and scenario:
	\begin{equation}
	g_\omega = l -  g^\mathrm{g}_\omega \; \forall \; \omega \in \Omega. \label{eq:power_to_residents}
	\end{equation}
	Note that this quantity is not necessarily the PV+S production since a portion of it may be directed to the grid as surplus.

	\subsection{CapEx and OpEx}
	
	We define the investor's CapEx as
	\begin{equation}
	\mathrm{CapEx} =  C^\mathrm{pv}_\mathrm{inv} + C^\mathrm{es}_\mathrm{inv} + \beta^\mathrm{pv} \cdot P^\mathrm{pv}_\mathrm{cap} + \beta^\mathrm{es} \cdot P^\mathrm{es}_\mathrm{cap} + C_\mathrm{grid}
	\end{equation}
	where $C^\mathrm{pv}_\mathrm{inv}$, $C^\mathrm{es}_\mathrm{inv}$, $P^\mathrm{pv}_\mathrm{cap}$, and $P^\mathrm{es}_\mathrm{cap}$ are introduced in Section~\ref{subsec:Sizing constraints}. The per-unit costs of PV and ES capacity are $\beta^\mathrm{pv}$ and $\beta^\mathrm{es}$, respectively. $C_\mathrm{grid}$ denotes a fixed grid connection cost. 
	
	We state the investor's OpEx as  
	\begin{equation}
	\mathrm{OpEx}_\omega=  \beta^\mathrm{es\; u}\cdot (c_\omega+d_\omega)+ \beta^\mathrm{mnt}\cdot P^\mathrm{pv}_\mathrm{cap}  + \tau^\mathrm{\intercal} g^\mathrm{s}_\omega\; \forall \; \omega \in \Omega \label{eq:def_OpEx}
	\end{equation}
	where the first term represents the ES utilization (i.e., degradation or ageing) cost as a function of charge and discharge; the second term represents the maintenance costs of the PV system; and the third term represents taxes for grid exports. The symbols $\beta^\mathrm{es\; u}$ and $\beta^\mathrm{mnt}$ in Eq.~\eqref{eq:def_OpEx} are the ES utilization coefficients (e.g., as in~\cite{Sarker_2017_Optimal}) and yearly maintenance cost per unit of installed PV, respectively.  The vector $\tau \in \mathbb R^T$ denotes taxes paid by the investor for exporting $g^\mathrm{s}_\omega$ to the grid.

	\section{Benefit allocation}
	\label{sec:benefit}
	
	The solution to the sizing problem is a Pareto-optimal PV+S system but does not tell us what the expected profit of the investor is nor what the expected savings of the consumers are. The reason is that the cost of PV+S energy for the consumer is part of the investor's revenue and cancels out in the objective.  Nevertheless, the solution to the sizing problem is important as it represents the maximum welfare PV+S system.
	
	From the sizing problem's solution, we define the \textit{net benefit}, the quantity to be shared among the investor and consumers, as follows.
	
	\begin{definition}
		The expected net benefit, $\overline B_\mathrm{net}$, is the difference between the present value of the consumers' electricity cost without PV+S, $C(l)$, and the objective of the sizing problem from Eq.~\eqref{eq:global_problem}, $\mathrm{obj}_\mathrm{LT}^*$.\label{def:net_benefit} 
		\begin{equation}
		\overline B_\mathrm{net} = C(l) - \mathrm{obj}_\mathrm{LT}^*
		\end{equation}
	\end{definition} 
	
	By definition, the expected net benefit is non-negative. This follows from the fact that not installing a PV+S system is a feasible solution (and therefore $C(l) \ge -\mathrm{obj}_\mathrm{LT}^*$). We can then say that any PV+S system will deliver positive expected benefits. The positivity of the expected net benefits is an important fact since it is a necessary condition for financial stability.
	
	The mechanisms whereby the net benefit is shared is the price of local energy $p$ and the amount of energy allocated to each consumer. A higher price $p$ translates to higher investor profit and lower savings for the consumers. Assuming that $p$ is lower than the average grid price, a consumer's savings are proportional to the amount of energy assigned to him/her.
	
	\subsection{The price of local energy}
	
	Let $\overline \Pi(p)$ denote expected investor profit as a function of the price $p$.  The relation 
	\begin{equation}
	\gamma \cdot \overline{B}_\mathrm{net} = \overline{\Pi}(p) \label{eq:NB_retained_by_investor}
	\end{equation}
	describes the share $\gamma$ of the net benefit kept by the investor as profit.  
	
	While it is relatively straight forward to calculate $\overline{B}_\mathrm{net}$ and $\overline{\Pi}(p)$ (by solving the sizing problem), calculating $\gamma$ or $p$ is not easy: there is no ``optimal'' point along the investor-consumer Pareto frontier. Therefore, there is no optimal $p$ nor a natural $\gamma$. 
	
	A possible scenario is that the investor and consumers arrive at a price through negotiation. Naturally, the negotiation does not happen inside a vacuum: it takes place in a particular economic, regulatory, technological, and cultural environment which makes the outcome hard to characterize for a general case. However, one can identify general dynamics. For example, we can argue that an abundance of capital diminishes the negotiating power of the investor. Conversely, access to effective PV production models by the consumer could allow them to better estimate project costs and benefits and claim a larger stake. While a deep study of the issue of benefit allocation is outside of this work, seminal work on bargaining, or negotiation, theory can be found in~\cite{nash1950bargaining, binmore1986nash} and an application to power systems in~\cite{Contreras_2019_Participation}.
	
	A basic characteristic of $\gamma$ is that it should normally be in $[0,1]$. An investor is unlikely to accept a negative profit, i.e., $\gamma<0$. Conversely, any $\gamma>1$ means that the consumers pay more than they would otherwise pay to the grid. In theory, $\gamma \notin [0,1]$ is possible if other (perhaps non-monetary) benefits are omitted from the sizing problem. For example, consumers could accept paying higher-than-grid prices for local and green energy. Conversely, an investor could accept losses if the project provides non-accounted benefits (e.g., a positive public image).

	\subsection{Energy allocation and key of repartition}
	\label{subsec:key}
	While $p$ determines the benefit allocation between investor and consumers \textit{as a whole}, the PV+S energy allocated to each consumer determines the benefit split \textit{among} the consumers themselves.  
	
	We allocate PV+S energy on three occasions throughout our framework. The first is on the planning stage before constituting the collective. There, we estimate the amount of local energy that each consumer can expect throughout the life of the project. With the estimate on hand, consumers can evaluate whether to enter the collective or not. We also allocate energy among consumers in the control and settlement stages. We detail these last two occasions in Section~\ref{sec:RTC}.
	
	What is important for consumers from a financial point of view is monthly or yearly energy allocation. However, under the French collective self-consumption framework, local energy production must be allocated among consumers in 30-minute intervals via a \textit{key of repartition}.

	In our work, a key of repartition has two fundamental characteristics. For every time interval: 
	\begin{enumerate}
		\item If there is an energy deficit (i.e., aggregate load $>$ PV+S production), all of the production must be assigned to the consumers and if there is a surplus, the aggregate load must be met by PV+S energy.
		\item The energy assigned to a consumer must be less than its load.
	\end{enumerate}    
	Let the matrix $G = [g_1, g_2,\hdots, g_N ]\in \mathbb R^{T\times N}_+$ denote a key of repartition\footnote{Recall that $g_i$ was introduced in Eq.~\eqref{eq:resident_cost} and is a vector energy directed to consumer $i$ during each time interval.}. The $(t,i)$\textsuperscript{th} element of $G$ represents the amount of energy that consumer $i$ receives during time interval $t$.  
	
	We formalize the two conditions above by defining the set
	\begin{equation}
	\mathcal G( g,L) = \left \{ G \mid G1 = \min\left( g,L1 \right), \; 0 \le G \le L \right \}, \label{eq:feasible_keys}
	\end{equation} where $g\in \mathbb R^T_+$ was introduced in Eq.~\eqref{eq:all_resident_cost} and is the total local energy assigned to consumers in each period. The matrix $L\in \mathbb R^{T\times N}_+$ is constructed by $L=[l_1, l_2, \hdots , l_N]$. Similar to $G$, the $(t,i)$\textsuperscript{th} element of $L$ is the load of consumer $i$ during interval $t$. The term $G1$  and $L1$ are $T$-long vectors composed of the sum of the columns of $G$ and $L$, respectively. The constraint\footnote{The operator $\min(\cdot, \cdot)$ is a element-wise function that takes the smallest number of each corresponding elements of the arguments. } $ G1 = \min\left( g,L1 \right)$ states Condition 1. The constraint $0\le G \le L$ states Condition 2.
	
	\begin{rmk}
		For any $g \in R_+^T$ and $L \in R_+^{T \times N}$ there is always a feasible key of repartition $G \in R_+^{T \times N}$. We analyze two cases: time periods with surplus PV generation ($g>L1$ and thus $L1=\mathrm{min}(g,L1)$) and time periods with deficit PV generation ($g\le L1$ which implies $g=\mathrm{min}(g,L1)$). For the first, the key of reparation $G=L$ satisfies both constraints of Eq.~\eqref{eq:feasible_keys}: $G1=L1=\mathrm{min}(g,L1)$ and $0\le G=L \le G$. For the second case, let $g=a \cdot L1$ where $0\le a < 1$. Then, we can build the key of repartition $G=a\cdot L$ which satisfies both constraints of Eq.~\eqref{eq:feasible_keys}: $G=a\cdot L=g=\mathrm{min}(g,L1)$ and $0\le G=a\cdot L \le L$. Since there are feasible keys of repartition for this two cases, we conclude that there is always a feasible key of repartition. 
	\end{rmk}
	
	For non-trivial cases, infinitely many keys of repartition satisfy $ \mathcal G$. In our work, we evaluate the set of keys on a \textit{basis of equity} and select the best. That is, we select the key from $\mathcal G$ that allocates yearly energy as equally as possible.

	Let the symbol $e \in \mathbb R_+^N$ 
	\begin{equation}
	e =  G^\intercal 1 \label{eq:def_e}
	\end{equation}
	denote yearly energy allocated to each consumer. The $i$\textsubscript{th} element of $e$ is the energy allocated to consumer $i$. Since equity in yearly allocation is our goal, the best key is one that assigns the same amount to each consumer, i.e., an $e$ whose elements have a zero variance.
	
	However, since the key is restricted by $\mathcal G$, it may not be possible to achieve a zero variance. Also, recall that $g$ is uncertain. Thus, we formulate a linearly-constrained quadratic optimization problem whose objective is to minimize the \textit{expected variance} of $e$.
	
	Let $\mathrm{Var}(e_\omega)$ denote the variance of the annual energy allotted to each consumer in scenario $\omega$. Based on work by Singh in~\cite{Singh_2019_Modeling} and Abada et al. in~\cite{abada2017viability}, we find the key that minimizes the expected variance of energy allocation by solving the convex\footnote{Let $\mu_e$ denote the mean value of the elements of $e$ and $\mathrm{Var}(e) =  \frac{1}{N} e^\intercal e - \mu_e^2 $ their variance. Since $\mathrm{Var}(e)$ is a quadratic function and $\mathcal E$ is defined by linear constraints, Problem~\eqref{eq:key_of_rep} is convex.} optimization problem
	\begin{equation}
	\min_{e_\omega \in \mathcal E(g_\omega, L_\omega) \forall \omega \in \Omega} \expected[\mathrm{Var}(e_{\omega})]. \label{eq:key_of_rep}
	\end{equation}
	Here, $ L_\omega = [l_{1,\omega}, \hdots,l_{N,\omega} ]^\intercal$ and $\mathcal E(g_\omega, L_\omega) = \{e \mid e=G^\intercal 1 , \;G \in \mathcal{G}(g_\omega,L_\omega) \}$.  The set $\mathcal E(g_\omega, L_\omega)$ constraints the annual energy allocation in each scenario, $e_\omega$, to be delivered by a feasible key of repartition. 
	
	The solution of Problem~\eqref{eq:key_of_rep}, $e_\omega^*$, serves to estimate the yearly amount of PV+S energy each consumer can expect during the life-time of the project: $\overline e = \expected[e_\omega^*]$. It also serves as a liaison between long-term planning stage and system operation: it is the reference value of how to allocate energy among consumers during operation.

	\section{System operation}
	\label{sec:RTC}

	The operation of the PV+S system has two objectives: to determine the control strategy and to formulate the \textit{actual} key of repartition. Whereas the key from Problem~\eqref{eq:key_of_rep} serves to estimate energy allocation, the actual key reflects system operation and the materialization of uncertainty. 
	
	We divide system operation into the control stage and the settlement stages. We formulate the control strategy on a rolling horizon fashion via an MPC algorithm. The goals of the control strategy are two-fold: to minimize the short-term costs of the collective and to formulate a key that tracks $\overline e$. In the settlement stage, we reconcile the MPC key with actual system operation and materialization of uncertainty by formulating yet another key, this time the definitive one. We formulate the definitive key via an optimization program that minimizes the mismatch between $\overline e$ and the sum of energy served and expected energy.

	\subsection{Control}
	\label{subsec:Control}
	Our MPC algorithm starts by determining a control strategy for a \textit{prediction horizon} (e.g., one or two days) by solving an optimization problem. Then, we implement the strategy during a shorter \textit{control horizon} (e.g., 30 minutes) and collect system operation data. Finally, we update predictions and the optimization problem and repeat the process.

	The core of MPC is the control problem solved at each iteration. It is very similar in structure to the sizing problem but differs in three main ways. First, while the sizing problem has a years-long horizon, the control problem's horizon is in the order of hours or days. Second, the control problem does not make sizing decisions since the PV+S system is already in place. Finally, the control problem has aims to track the allocation of energy $\overline e$.  
	
	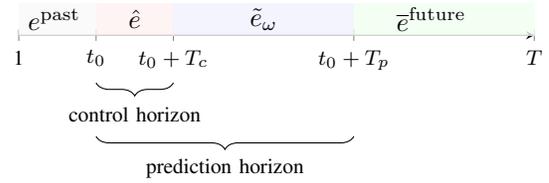
\begin{figure}
\begin{center}
\begin{tikzpicture}
\begin{axis}[%
axis x line=center,
axis y line=none,
xmin=0,xmax=10, ymin=-3, ymax=1,
xtick={0,1.5,3,6.5,10},
xticklabels={1,$t_0$,$t_0+T_c$,$t_0+T_p$,$T$},
]
\fill[gray!5, opacity=0.8] (0,0) rectangle (1.5,0.3);
\fill[red!5, opacity=0.8] (1.5,0) rectangle (3,0.3);
\fill[blue!5, opacity=0.8] (3,0) rectangle (6.5,0.3);
\fill[green!5, opacity=0.8] (6.5,0) rectangle (10,0.3);
\node[align=left, anchor=west] at (0,.15) {$e^\mathrm{past}$};
\node[align=left, anchor=center] at (2.25,.15) {$\hat e$};
\node[align=left, anchor=center] at (4.75,.15) {$\tilde e_\omega$};
\node[align=left, anchor=center] at (8,.15) {$\overline e^{\mathrm{future}}$};
\draw [decorate,decoration={brace,amplitude=5pt,mirror,raise=4ex}]
  (1.5,0) -- (3,0) node[midway,yshift=-3em]{{\footnotesize control  horizon}};
  \draw [decorate,decoration={brace,amplitude=5pt,mirror,raise=4ex}]
  (1.5,-.5) -- (6.5,-.5) node[midway,yshift=-3em]{{\footnotesize prediction  horizon}};
\end{axis}

\end{tikzpicture}
\end{center}
\caption{Illustration of the control and prediction horizon. } \label{fig:control_prediction_horizons}
\end{figure}

	As illustrated in Fig.~\ref{fig:control_prediction_horizons}, let $[t_0,t_0+T_c]$ and $[t_0,t_0+T_p]$ denote the prediction and control horizons, respectively.  We model decisions in the control horizon as first stage variables and identify them with a \textit{hat}. Decisions in the rest of the prediction horizon, i.e., in $[t_0+T_c,t_0+T_p]$, are second stage variables and are identified by a \textit{tilde}. For example, PV production during the control horizon, $[t_0,t_0+T_c]$ is denoted by $\hat g^\mathrm{pv} \in \mathbb R_+^{T_p}$ and during $[t_0+T_c,t_0+T_p]$ as  $\tilde g_\omega^\mathrm{pv} \in \mathbb R_+^{T_p}$. We represent uncertainty in the control horizon with a central forecast and via scenarios in the rest of the prediction horizon. 
	
	We write the control problem as 
	\begin{subequations}
		\begin{align}
		\min &  \expected\!\left[ C\!\left(\!\begin{bmatrix} \hat l \\ \tilde l_\omega \end{bmatrix}\! -\!  \begin{bmatrix} \hat g \\ \tilde g_\omega \end{bmatrix}\!\right) \!+\! \mathrm{OpEx}_{\omega} \!-\! \begin{bmatrix} \hat \lambda \\ \tilde \lambda_\omega \end{bmatrix}^\intercal\! \begin{bmatrix} \hat g^\mathrm{s}  \\ \tilde g^\mathrm{s} _\omega \end{bmatrix}\right] \!+\!   \theta\! \cdot\!  \|  \overline m\|_2^2  \nonumber \\
		\mathrm{s.t.} \;      & \mathrm{Eqs.}~\eqref{eq:def_g_pv},~\eqref{eq:def_ES_const},~\eqref{eq:energy_direction_def},~\eqref{eq:power_to_residents},~\eqref{eq:def_OpEx},~\eqref{eq:def_e} \\ 
		&  \hat e + \tilde e_\omega \in \mathcal E \left(\begin{bmatrix} \hat g \\ \tilde g_\omega \end{bmatrix}, \begin{bmatrix}\hat L \\ \tilde L_\omega \end{bmatrix}\right)  \; \forall  \; \omega \in \Omega \label{eq:def_e_tilde}\\
		& \overline m = \expected\left[e^{\mathrm{past}}+ \hat e+ \tilde e_\omega +  \overline e^{\mathrm{future}} - \overline{e} \right]. \ \label{eq:def_mismatch}
		\end{align} \label{eq:mpc_prob}
	\end{subequations}
	The objective has two components. The first is the operation costs during the control horizon (denoted by symbols with a hat) and during the rest of the prediction horizon (denoted by symbols with a tilde). $\mathrm{OpEx}_\omega$ denotes PV+S OpEx during the entire prediction horizon. The second is a term proportional to the $\ell^2$ norm of the expected mismatch, $\overline m$, between $\overline e$ from long-term planning and actual PV+S energy delivered to consumers. The mismatch is weighted by a parameter $\theta \in \mathbb R_+$ that has units of $\frac{\text{EUR}}{\text{kWh}^2}$. The paragraph that follow describes $\overline m$ in greater detail.  
	
	The solution space is subject to operation constraints, energy flow definitions, and OpEx definition defined by Eqs.~\eqref{eq:def_g_pv}-\eqref{eq:def_e}. Constraints~\eqref{eq:def_e_tilde} constraints energy allocated during the prediction horizon to be in $\mathcal E$. Finally, Eq.~\eqref{eq:def_mismatch} defines the expected mismatch as the difference between the promised energy to consumers $\overline e$ and
	\begin{itemize}
		\item $e^{\mathrm{past}}$, the energy allocated to each consumer prior to $t_0$;
		\item $\hat e$, allocated energy during the control horizon;
		\item $\tilde e_\omega$, allocated energy during $[t_0 + T_c, t_0+T_p]$;
		\item $\overline e^{\mathrm{future}}_\omega$, the expected\footnote{In our work, this expectation is obtained from $G^*_\omega$ from Problem~\eqref{eq:key_of_rep}} energy for each consumer after $T_p$.
	\end{itemize}

	\subsection{Settlement}
	\label{sec:PS}
	
	Uncertainty during the control horizon is relatively small since we solve the control problem shortly before $t_0$ when good predictions are available. Nevertheless, uncertain parameters are still present and materialize sometime before the end of the control horizon. Thus, we formulate an ex-post operation settlement problem to reconcile the key of repartition from the control strategy and actual operation of the PV+S system. 
	
	Let $\hat g^*$ denote energy to be assigned to consumers by Problem~\eqref{eq:mpc_prob} and $\hat \epsilon$ denote the deviation between actual and expected generation. Then, the energy assigned to consumers is $[\hat g^* + \hat \epsilon]^+$. We take the positive part of $\hat g^* + \hat \epsilon$ because no ``negative" energy can be assigned to consumers. Negative elements of $\hat g^* + \hat \epsilon$ is energy bought from the grid.   
	
	We formulate the settlement problem as
	\begin{equation*} \min_{\hat e \in \mathcal E ([\hat g^* + \hat \epsilon]^+, \hat L)}     \| \expected[e^{\mathrm{past}}+ \hat e + \tilde e_\omega +  \overline e^{\mathrm{future}} - \overline{e} ]\|_2^2.\end{equation*}
	Similar to the control problem, its objective is to minimize the $\ell^2$ norm of expected mismatch between promised energy to consumers, $\overline e$, and the energy assigned to each one of them during the four divisions of time illustrated in Fig.~\ref{fig:control_prediction_horizons}. In this case, $\hat e$ is the only optimization variable and that the objective is to allocate $[\hat g^* + \hat \epsilon]^+$.
	
	\section{Case Study}
	\label{sec:casestudy}
	\subsection{Background and data}
	
	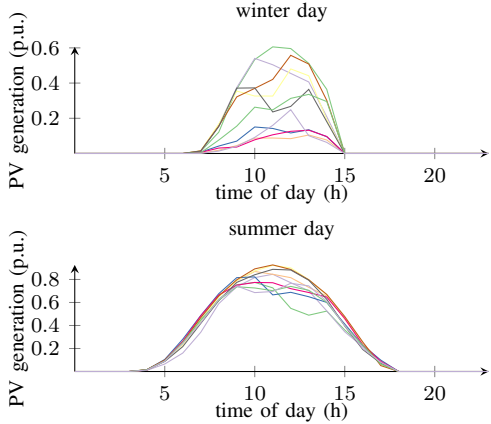
\begin{figure}
\centering
\begin{tikzpicture}

    \begin{axis}[title=\footnotesize{winter day},
    cycle list/Accent,
            xmax=23,
            legend style={at={(0.9,0.7)},
         anchor=north,legend columns=1,draw=none, font=\footnotesize},
            ticklabel style = {font=\footnotesize},
            width=0.8\linewidth, height=3cm,    
          ylabel= PV generation (p.u.),
            xlabel= time of day (h),
              y label style={at={(axis description cs:-0.09,.5)},rotate=90,anchor=south,align=center, font=\footnotesize},
              x label style={at={(axis description cs:0.5,-0.2)},anchor=north, font=\footnotesize},
         ]        

\addplot table[x = time, y=s1] {solar_scenarios.dat};
\addplot table[x = time, y=s2] {solar_scenarios.dat};
\addplot table[x = time, y=s3] {solar_scenarios.dat};
\addplot table[x = time, y=s4] {solar_scenarios.dat};
\addplot table[x = time, y=s5] {solar_scenarios.dat};
\addplot table[x = time, y=s6] {solar_scenarios.dat};
\addplot table[x = time, y=s7] {solar_scenarios.dat};
\addplot table[x = time, y=s8] {solar_scenarios.dat};
\addplot table[x = time, y=s9] {solar_scenarios.dat};
\addplot table[x = time, y=s10] {solar_scenarios.dat};

\end{axis}

\end{tikzpicture}

\begin{tikzpicture}

    \begin{axis}[title=\footnotesize{summer day},
    cycle list/Accent,
            xmax=23,
            legend style={at={(0.9,0.7)},
         anchor=north,legend columns=1,draw=none, font=\footnotesize},
            ticklabel style = {font=\footnotesize},
            width=0.8\linewidth, height=3cm,    
          ylabel= PV generation (p.u.),
            xlabel= time of day (h),
              y label style={at={(axis description cs:-0.09,.5)},rotate=90,anchor=south,align=center, font=\footnotesize},
              x label style={at={(axis description cs:0.5,-0.2)},anchor=north, font=\footnotesize},
         ]        

\addplot table[x = time, y=s1] {solar_scenarios_summer.dat};
\addplot table[x = time, y=s2] {solar_scenarios_summer.dat};
\addplot table[x = time, y=s3] {solar_scenarios_summer.dat};
\addplot table[x = time, y=s4] {solar_scenarios_summer.dat};
\addplot table[x = time, y=s5] {solar_scenarios_summer.dat};
\addplot table[x = time, y=s6] {solar_scenarios_summer.dat};
\addplot table[x = time, y=s7] {solar_scenarios_summer.dat};
\addplot table[x = time, y=s8] {solar_scenarios_summer.dat};
\addplot table[x = time, y=s9] {solar_scenarios_summer.dat};
\addplot table[x = time, y=s10] {solar_scenarios_summer.dat};

\end{axis}

\end{tikzpicture}

  \caption{The upper plot shows PV generation scenarios (in per unit of PV capacity) for a winter day. The lower plot shows scenarios for a summer day.  }
\label{fig:solar_data}
\end{figure}
	
	We demonstrate our proposed framework with a case study of a collective of 15 residential consumers sharing a PV+S system. We use 10 scenarios of solar irradiation for the south of France from~\cite{jensen_RE-Europe_data}. Data for two sample days (one winter one summer) is shown in Fig.~\ref{fig:solar_data}. French electric system fees and taxes (\textit{TURPE}) are  from~\cite{deliberation_CRE} and energy prices from \textit{\'Electrcit\'e de France's (EDF) Tarif Bleu}\footnote{The \textit{tarif bleu} includes the taxes and fees so we isolate the energy price by subtracting them.}. We use load data of residences in San Diego, California from~\cite{pecan_street}.

	We consider solar PV and lithium-ion ES investments and analyze two cases: a baseline case and a pessimistic case. In the former, we consider low per-kW of PV prices\footnote{We use an exchange rate of \$1=1.11 EUR.}  (1.1 EUR/W for systems under 100 kW and 0.95 EUR/W for PV systems over 100 kW~\cite{photovoltaique.com}) and a price of 0.06 EUR/kWh for surplus sold energy to the grid~\cite{edf_france_2018}. In the later, we assume high per-kW PV prices (1.68 EUR/W for systems  $<$100 KW and 1.58 EUR/W for PV systems $\ge$100 kW~\cite{fu2018us}) and no compensation for grid injections. 
	
	We consider 4 available inverters with capacities of 50, 99, 157, and 249 kW and costs of 72 EUR/kW~\cite{fu2018us}. We assume a per-kWh cost of ES of 158EUR~\cite{Lavergne_2019_Stockage} and a round-trip efficiency of 90\%. We consider PV subsidies of 100 EUR/kW for systems under 100 kW~\cite{edf_france_2018}. The planning horizon is 20 years and is divided into 1-hour periods in the sizing problem and into 30-minute periods in the control problem. We assume a discount rate of 3\%.
	
	We ran the simulations on a Windows 10 laptop computer running on an Intel\textcopyright~Core\texttrademark~i5-8250U CPU @1.60GHz 1.80GHz with 8GB of RAM. The models are coded in Julia and the optimization problems solved using Gurobi 8.0.1.

	\subsection{Sizing}
	In both the baseline and pessimistic cases, the sizing problem contains 526,152 variables (10 of which are binary) and 701,337 constraints. The binary variables describe discrete choices (e.g., inverter) and non-convexities\footnote{The PV cost and subsidy functions are non-convex.}. Despite its size and non-convexities, the problem is relatively easy to solve: the baseline case solves in approximately 2 minutes and the pessimistic in 12 minutes. 
	
	\begin{figure}
\begin{center}
\begin{tikzpicture}
    \begin{axis}[
            xmax=100,
            legend style={at={(0.9,0.7)},
         anchor=north,legend columns=1,draw=none, font=\footnotesize},
            ticklabel style = {font=\footnotesize},
            width=0.8\linewidth, height=3.5cm,    
          ylabel= expected cost \\ (\faEuro/month/consumer),
            xlabel= PV capacity (kW),
              y label style={at={(axis description cs:-0.09,.5)},rotate=90,anchor=south,align=center, font=\footnotesize},
              x label style={at={(axis description cs:0.5,-0.15)},anchor=north, font=\footnotesize},
         ]        
        \addplot[color=red, thick] table[x=x,y=pes]  {long_term_plan.dat};
        \addplot[color=blue, thick] table[x=x,y=opt]  {long_term_plan.dat};
        
        \legend{pessimistic case, baseline case}
        
    \end{axis}
\end{tikzpicture}
\end{center}
    \caption{Expected monthly cost per consumer (in present value) as a function of PV capacity. The optimum for the pessimistic case is 50 kW and 249 kW for the baseline case (not shown in the plot). }

  \label{fig:long_term_cost}
\end{figure}
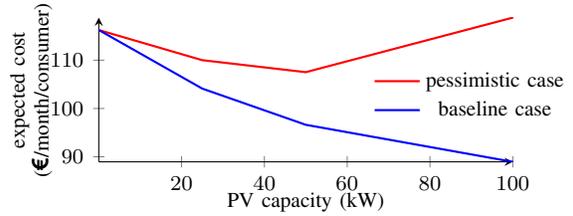
	\begin{table}[ht]
\caption{Expected costs, benefit, and net benefit in both cases.}
\begin{center}
\begin{tabular}{|c|c|c|c|c|c|}
\hline
\multicolumn{2}{|c|}{\textbf{costs (k\faEuro)}}               & \multicolumn{3}{c|}{\textbf{benefit (k\faEuro)}}                                          &                                                                                       \\ \cline{1-5}
\textbf{CapEx}             & \textbf{OpEx}             & \textbf{grid savings}      & \textbf{subsidy}         & \textbf{surplus}           & \multirow{-2}{*}{\textbf{\begin{tabular}[c]{@{}c@{}}net benefit \\ (k\faEuro)\end{tabular}}} \\ \hline
\multicolumn{6}{|c|}{\textbf{baseline case}}                                                                                                                                                                                      \\ \hline
{\color[HTML]{CB0000} 258} & {\color[HTML]{CB0000} 86} & {\color[HTML]{009901} 224} & 0                        & {\color[HTML]{009901} 262} & {\color[HTML]{009901} 142}                                                            \\ \hline
\multicolumn{6}{|c|}{\textbf{pessimistic case}}                                                                                                                                                                                     \\ \hline
{\color[HTML]{CB0000} 110} & {\color[HTML]{CB0000} 10} & {\color[HTML]{009901} 148} & {\color[HTML]{009901} 5} & 0                          & {\color[HTML]{009901} 32}                                                             \\ \hline
\end{tabular}
\end{center}
\label{table:cost_benefits}
\end{table}

	Fig.~\ref{fig:long_term_cost} shows the expected consumer monthly cost (in present value) as function of PV capacity. In the baseline case, the optimal PV capacity is 249 kW and the optimal ES capacity is 150 kWh. In this case, we expect the PV+S system to lower costs from 116 EUR/consumer/month to 77 EUR/consumer/month. This represents net benefits of 142 thousand EUR~during the lifetime of the project. The pessimistic case (50kW PV/127 kWh ES), on the other hand, represents lifetime net benefits of 32 thousand EUR. Table~\ref{table:cost_benefits} breaks down costs and benefits for both cases. Hereafter we analyze the pessimistic case.

	\subsection{Benefit allocation}
	
	\begin{table}[ht]
    \caption{Investor profit and consumer savings under different prices for both the baseline and the optimistic cases.}
\textbf{\underline{Pessimistic case}} 
  \begin{center}
\begin{tabular}{c|c|c|c} 
\textbf{$\boldsymbol{p}$ (\faEuro/kWh)} & $\boldsymbol{\gamma}$ & \textbf{\begin{tabular}[c]{@{}c@{}} investor profit (k\faEuro) \end{tabular}} & \multicolumn{1}{c}{\textbf{\begin{tabular}[c]{@{}c@{}}  consumer  savings (k\faEuro) \end{tabular}}} \\
\hline
0.10             & 0              & 0                                                                               & {\color[HTML]{009901}32 }                                                                                                           \\ 
0.115             & 0.5            & {\color[HTML]{009901}16 }                                                                             & {\color[HTML]{009901}16   }                                                                                                           \\ 
0.13             & 1              & {\color[HTML]{009901}32  }                                                                           & 0                                                                                                               \\ 
\end{tabular}

\end{center}
\textbf{\underline{Baseline case}} 
  \begin{center}
\begin{tabular}{c|c|c|c} 
\textbf{$\boldsymbol{p}$ (\faEuro/kWh)} & $\boldsymbol{\gamma}$ & \textbf{\begin{tabular}[c]{@{}c@{}} investor profit (k\faEuro) \end{tabular}} & \multicolumn{1}{c}{\textbf{\begin{tabular}[c]{@{}c@{}}  consumer  savings (k\faEuro) \end{tabular}}} \\
\hline
0.05             & 0              & 0                                                                               & {\color[HTML]{009901}142 }                                                                                                           \\ 
0.09             & 0.5            & {\color[HTML]{009901}71 }                                                                             & {\color[HTML]{009901}71   }                                                                                                           \\ 
0.13             & 1              & {\color[HTML]{009901}142  }                                                                           & 0                                                                                                               \\ 
\end{tabular}
\end{center}

    \label{table:benefit_split}
    
\end{table}

	As discussed in Section~\ref{sec:benefit}, the price of PV+S energy $p$ is the main mechanism used to allocate benefits between the investor and the consumers. As shown in Table~\ref{table:benefit_split}, the ``break-even'' price for the investor is 0.10 EUR/kWh in the baseline case and 0.05 EUR/kWh in the optimistic case. On the other hand, the price at which consumers can expect to pay the same whether or not they join the collective is 0.13 EUR/kWh for both cases. Thus, any price between the investor break-even and the consumer break-even prices allocates \textit{positive benefits} to both the investor and the consumers (i.e., a win-win situation). We call this price range the win-win price range. 
	
	Note that while the consumer break-even price is the same in both cases, for the investor it is higher in the baseline case. The reason the consumers' break-even price is constants is that their grid costs (the alternative to energy from the PV+S system) are the same in both cases. For the investor, however, the higher CapEx and OpEx in the baseline case (see Table~\ref{table:cost_benefits}) increase the break-even price. In general, higher capital and operational costs lower the expected net benefit and shrink the win-win price range. On the other hand, more favorable conditions for PV+S investment (e.g., higher grid prices, lower CapEx and OpEx) expand the win-win price range. 
	
	While predicting $p$ is outside the scope of this work, identifying the win-win price range is important because it defines the financially sustainable range of $p$. Hereafter we analyze the pessimistic case where benefits are allocated equally among the investor and consumers, i.e., $p=0.115$.

	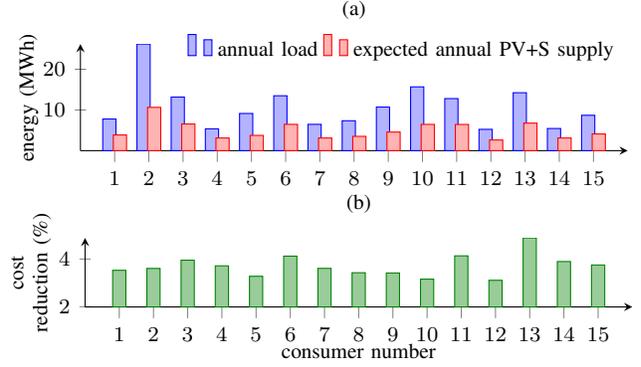
\begin{figure}
\centering
 \begin{tikzpicture}
\begin{axis}[
title = \footnotesize{(a)},
    legend cell align={left},
    ymin=0,  xmin=0, xmax=16,
    width=\linewidth, height=3cm, 
	ticklabel style = {font=\footnotesize},
	y label style={at={(axis description cs:-0.06,.5)},rotate=90,anchor=south,align=center, font=\footnotesize},
    x label style={at={(axis description cs:0.5,-0.2)},anchor=north, font=\footnotesize},
	ylabel=energy (MWh), 
	legend style={at={(1,0.95)}, draw=none,fill=none,
		anchor=east,legend columns=-1, font=\footnotesize},
	ybar=-1pt,
	bar width=5pt,
	xtick={1,2,3,4,5,6,7,8,9,10,11,12,13,14,15},
]

\addplot table[x=x, y=annual_load] {savings_per_cons.dat};
\addplot table[x=x, y=overline_e] {savings_per_cons.dat};

\legend{annual load , expected annual PV+S supply}

\end{axis}

\end{tikzpicture}
\vspace{-3pt}

 \begin{tikzpicture}

\definecolor{ao(english)}{rgb}{0.0, 0.5, 0.0}
\begin{axis}[
title = \footnotesize{(b)},
    legend cell align={left},
    ymin=1.99,  xmin=0, xmax=16,
    width=\linewidth, height=2.5cm, 
	ticklabel style = {font=\footnotesize},
	y label style={at={(axis description cs:-0.045,.5)},rotate=90,anchor=south,align=center, font=\footnotesize},
    x label style={at={(axis description cs:0.5,-0.4)},anchor=north, font=\footnotesize},
	ylabel=cost \\ reduction (\%), xlabel=consumer number ,
	ybar=-2pt,
	bar width=5pt,
	xtick={1,2,3,4,5,6,7,8,9,10,11,12,13,14,15},
	ytick={2,4},
]

\addplot[ao(english), fill=ao(english)!40] table[x=x, y=cost_reduction] {savings_per_cons.dat};
\end{axis}
\end{tikzpicture}
  \caption{Plot (a) shows annual load (blue) and expected PV+S energy (red) for each of the 15 residents. Plot (b) shows the expected percentage cost reduction from entering the collective with $p=0.115$.}

\label{fig:benefit_allocation}
\end{figure}
	
	Each consumer's savings is a function of the amount of PV+S energy each receives at a price $p$. As discussed in Section~\ref{sec:benefit}, infinitely many energy allocations are possible but we use Problem~\eqref{eq:key_of_rep} to find a key of repartition (i.e., the PV+S supply allotment to each consumer on a 30-minute basis) that allocates energy equitably. The red bars in Fig.~\ref{fig:benefit_allocation}(a) show the annual PV+S supply assigned to each consumer by the solution of Problem~\eqref{eq:key_of_rep}.
	
	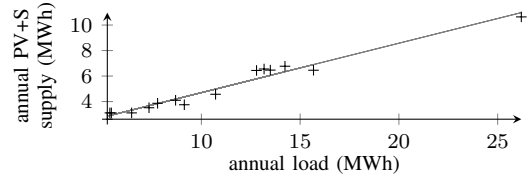
\begin{figure}
\centering
\begin{tikzpicture}
\begin{axis}[          ylabel= annual PV+S \\ supply (MWh),
            xlabel= annual load (MWh),
                          y label style={at={(axis description cs:-0.1,.5)},rotate=90,anchor=south,align=center, font=\footnotesize},
              x label style={at={(axis description cs:0.5,-0.25)},anchor=north, font=\footnotesize},
              width=0.8\linewidth, height=3cm,
              legend style={at={(1,0.3)},
         anchor=east,legend columns=1,draw=none, font=\footnotesize},
]
\addplot[mark=none, gray] table[x=annual_load, y={create col/linear regression={y=overline_e}}] {savings_per_cons.dat};
\addplot [only marks, mark = +] table[x= annual_load, y=overline_e] {savings_per_cons.dat};
\end{axis}
\end{tikzpicture}
\caption{Expected annual load and PV+S supply. Each data-point represents a consumer and the line is the least-square error linear approximation (R-squared = 0.97).}
\label{fig:regression}
\end{figure}

	Even though the objective of Problem~\eqref{eq:key_of_rep} is to equally distribute PV+S supply, one can notice that the supply distribution is not quite equal. The reason is that the restrictions imposed by the set $\mathcal G$, which defines a key of repartition, prevents the allotment from being perfectly equitable. For example, a load that systematically ``under-consumes'' during high PV production hours tends to ``miss out'' on the opportunity to receive PV+S energy. In fact, by examining the blue and red bars of~\ref{fig:benefit_allocation}(a) one can see that the expected PV+S supply is closely linked to the annual load. Furthermore, as shown by Fig.~\ref{fig:regression} most of the variance in annual PV+S supply can be explained by the consumer's annual load (with an R-squared = 0.97). What this means is that, as shown in Fig.~\ref{fig:benefit_allocation}(b) cost reductions in terms of percentage are relatively uniform among consumers and most consumers can expect a 3 to 5\% cost reduction by joining the collective.

	\subsection{Operation and settlement}
	\begin{figure}
\centering
\begin{tikzpicture}

\begin{axis}[
            legend style={draw=none, font=\footnotesize},
cycle list/Accent,
  date coordinates in=x,
  ticklabel style = {font=\footnotesize},
  xticklabel=\day-\month,
  xmin  = 00-01-01 01:00, xmax = 00-12-30 19:00,
  ymin=-15,
  width=.95\linewidth, height=3cm, 
  y label style={at={(axis description cs:-0.1,.5)},rotate=90,anchor=south,align=center, font=\footnotesize},
    x label style={at={(axis description cs:0.5,-0.2)},anchor=north, font=\footnotesize},
    axis x line*=bottom,
    ylabel=mismatch (\%), xlabel= day-month,
    xtick = {00-03-01 01:00, 00-06-01 01:00, 00-09-01 01:00, 00-012-01 01:00},
    title = \footnotesize{(a) proposed algorithm}
]

\addplot[mark=none, gray, dashed] coordinates {(00-01-00 01:00,6.715883
) (00-12-30 01:00,6.715883
)};
\addplot table[x = time, y expr=\thisrowno{1}*100, col sep=comma] {energy_mismatch_csv.dat};
\addplot table[x = time, y expr=\thisrowno{2}*100, col sep=comma] {energy_mismatch_csv.dat};
\addplot table[x = time, y expr=\thisrowno{3}*100, col sep=comma] {energy_mismatch_csv.dat};
\addplot table[x = time, y expr=\thisrowno{4}*100, col sep=comma] {energy_mismatch_csv.dat};
\addplot table[x = time, y expr=\thisrowno{5}*100, col sep=comma] {energy_mismatch_csv.dat};
\addplot table[x = time, y expr=\thisrowno{6}*100, col sep=comma] {energy_mismatch_csv.dat};
\addplot table[x = time, y expr=\thisrowno{7}*100, col sep=comma] {energy_mismatch_csv.dat};
\addplot table[x = time, y expr=\thisrowno{8}*100, col sep=comma] {energy_mismatch_csv.dat};
\addplot table[x = time, y expr=\thisrowno{9}*100, col sep=comma] {energy_mismatch_csv.dat};
\addplot table[x = time, y expr=\thisrowno{10}*100, col sep=comma] {energy_mismatch_csv.dat};
\addplot table[x = time, y expr=\thisrowno{11}*100, col sep=comma] {energy_mismatch_csv.dat};
\addplot table[x = time, y expr=\thisrowno{12}*100, col sep=comma] {energy_mismatch_csv.dat};
\addplot table[x = time, y expr=\thisrowno{13}*100, col sep=comma] {energy_mismatch_csv.dat};
\addplot table[x = time, y expr=\thisrowno{14}*100, col sep=comma] {energy_mismatch_csv.dat};
\addplot table[x = time, y expr=\thisrowno{15}*100, col sep=comma] {energy_mismatch_csv.dat};

\legend{PV generation mismatch}

\end{axis}

\end{tikzpicture}

\begin{tikzpicture}

\begin{axis}[
            legend style={draw=none, font=\footnotesize},
cycle list/Accent,
  date coordinates in=x,
  ticklabel style = {font=\footnotesize},
  xticklabel=\day-\month,
  xmin  = 00-01-01 01:00, xmax = 00-12-30 19:00,
  ymin=-15,
  width=.95\linewidth, height=3cm, 
  y label style={at={(axis description cs:-0.1,.5)},rotate=90,anchor=south,align=center, font=\footnotesize},
    x label style={at={(axis description cs:0.5,-0.2)},anchor=north, font=\footnotesize},
    axis x line*=bottom,
    ylabel=mismatch (\%), xlabel= day-month,
    xtick = {00-03-01 01:00, 00-06-01 01:00, 00-09-01 01:00, 00-012-01 01:00},
    title=\footnotesize{ (b) MPC, myopic (short-sighted) settlement}
]

\addplot[mark=none, gray, dashed] coordinates {(00-01-00 01:00,6.715883
) (00-12-30 01:00,6.715883
)};
\addplot table[x = time, y expr=\thisrowno{1}*100, col sep=comma] {energy_mismatch_myopic.dat};
\addplot table[x = time, y expr=\thisrowno{2}*100, col sep=comma] {energy_mismatch_myopic.dat};
\addplot table[x = time, y expr=\thisrowno{3}*100, col sep=comma] {energy_mismatch_myopic.dat};
\addplot table[x = time, y expr=\thisrowno{4}*100, col sep=comma] {energy_mismatch_myopic.dat};
\addplot table[x = time, y expr=\thisrowno{5}*100, col sep=comma] {energy_mismatch_myopic.dat};
\addplot table[x = time, y expr=\thisrowno{6}*100, col sep=comma] {energy_mismatch_myopic.dat};
\addplot table[x = time, y expr=\thisrowno{7}*100, col sep=comma] {energy_mismatch_myopic.dat};
\addplot table[x = time, y expr=\thisrowno{8}*100, col sep=comma] {energy_mismatch_myopic.dat};
\addplot table[x = time, y expr=\thisrowno{9}*100, col sep=comma] {energy_mismatch_myopic.dat};
\addplot table[x = time, y expr=\thisrowno{10}*100, col sep=comma] {energy_mismatch_myopic.dat};
\addplot table[x = time, y expr=\thisrowno{11}*100, col sep=comma] {energy_mismatch_myopic.dat};
\addplot table[x = time, y expr=\thisrowno{12}*100, col sep=comma] {energy_mismatch_myopic.dat};
\addplot table[x = time, y expr=\thisrowno{13}*100, col sep=comma] {energy_mismatch_myopic.dat};
\addplot table[x = time, y expr=\thisrowno{14}*100, col sep=comma] {energy_mismatch_myopic.dat};
\addplot table[x = time, y expr=\thisrowno{15}*100, col sep=comma] {energy_mismatch_myopic.dat};


\end{axis}

\end{tikzpicture}

\begin{tikzpicture}

\begin{axis}[
            legend style={draw=none, font=\footnotesize},
cycle list/Accent,
  date coordinates in=x,
  ticklabel style = {font=\footnotesize},
  xticklabel=\day-\month,
  xmin  = 00-01-01 01:00, xmax = 00-12-30 19:00,
  ymin=-15,
  width=.95\linewidth, height=3cm, 
  y label style={at={(axis description cs:-0.1,.5)},rotate=90,anchor=south,align=center, font=\footnotesize},
    x label style={at={(axis description cs:0.5,-0.2)},anchor=north, font=\footnotesize},
    axis x line*=bottom,
    ylabel=mismatch (\%), xlabel= day-month,
    xtick = {00-03-01 01:00, 00-06-01 01:00, 00-09-01 01:00, 00-012-01 01:00},
    title=\footnotesize{ (c) rule-based control, myopic settlement}
]

\addplot[mark=none, gray, dashed] coordinates {(00-01-00 01:00,6.715883
) (00-12-30 01:00,6.715883
)};
\addplot table[x = time, y expr=\thisrowno{1}*100, col sep=comma] {energy_mismatch_myopic_rule-based.dat};
\addplot table[x = time, y expr=\thisrowno{2}*100, col sep=comma] {energy_mismatch_myopic_rule-based.dat};
\addplot table[x = time, y expr=\thisrowno{3}*100, col sep=comma] {energy_mismatch_myopic_rule-based.dat};
\addplot table[x = time, y expr=\thisrowno{4}*100, col sep=comma] {energy_mismatch_myopic_rule-based.dat};
\addplot table[x = time, y expr=\thisrowno{5}*100, col sep=comma] {energy_mismatch_myopic_rule-based.dat};
\addplot table[x = time, y expr=\thisrowno{6}*100, col sep=comma] {energy_mismatch_myopic_rule-based.dat};
\addplot table[x = time, y expr=\thisrowno{7}*100, col sep=comma] {energy_mismatch_myopic_rule-based.dat};
\addplot table[x = time, y expr=\thisrowno{8}*100, col sep=comma] {energy_mismatch_myopic_rule-based.dat};
\addplot table[x = time, y expr=\thisrowno{9}*100, col sep=comma] {energy_mismatch_myopic_rule-based.dat};
\addplot table[x = time, y expr=\thisrowno{10}*100, col sep=comma] {energy_mismatch_myopic_rule-based.dat};
\addplot table[x = time, y expr=\thisrowno{11}*100, col sep=comma] {energy_mismatch_myopic_rule-based.dat};
\addplot table[x = time, y expr=\thisrowno{12}*100, col sep=comma] {energy_mismatch_myopic_rule-based.dat};
\addplot table[x = time, y expr=\thisrowno{13}*100, col sep=comma] {energy_mismatch_myopic_rule-based.dat};
\addplot table[x = time, y expr=\thisrowno{14}*100, col sep=comma] {energy_mismatch_myopic_rule-based.dat};
\addplot table[x = time, y expr=\thisrowno{15}*100, col sep=comma] {energy_mismatch_myopic_rule-based.dat};


\end{axis}

\end{tikzpicture}

  \caption{Cumulative mismatch as a function of time for each consumer under three algorithms: (a) the proposed framework, (b) the proposed MPC and a myopic (short-sighted) settlement, and (c) a rule-based control algorithm and the myopic settlement. The proposed framework delivers smaller end-of year mismatches. }
 \label{fig:energy_mismatch}
\end{figure}
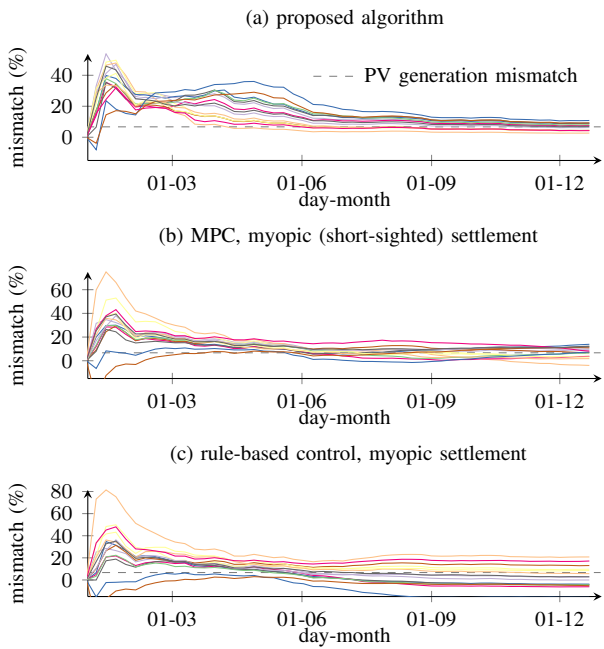
	
	The objective of the operation phase is to find the control strategy that minimizes short-term costs and mismatch between promised energy (the red bars of Fig.~\ref{fig:benefit_allocation}(a)) and actual delivered energy to each consumer.
	
	Fig.~\ref{fig:energy_mismatch} shows the consumer mismatch in terms of percentage energy promised as a function of time for three different algorithms. Plot (a) shows results using the proposed MPC and settlement algorithm, plot (b) shows results using the proposed MPC algorithm and a myopic (short-sighted) settlement algorithm, and plot (c) shows results using a rule-based control algorithm similar to the one in~\cite{weniger_2014_sizing} and the myopic settlement algorithm. The myopic settlement allocates energy to the consumers as fairly as possible but does not take into consideration past nor expected future allocations. Under the rule-based control, the ES system operates under the following regime: when there is a PV energy surplus, it charges as much as possible (while observing power and energy limits) and when there is a PV energy deficit the ES discharges as much as possible, also while observing limits. Note that the most important flaw of the rule-based algorithm is that it does not seek to minimize costs as it fails to account for retail energy prices. 
	
	\begin{figure}
\centering

 \begin{tikzpicture}

\definecolor{ao(english)}{rgb}{0.0, 0.5, 0.0}
\begin{axis}[
    legend cell align={left},
    width=.85\linewidth, height=2.5cm, 
	ticklabel style = {font=\footnotesize},
	y label style={at={(axis description cs:-0.1,.5)},rotate=90,anchor=south,align=center, font=\footnotesize},
    x label style={at={(axis description cs:0.5,-0.3)},anchor=north, font=\footnotesize},
	ylabel=cumulative \\end-of-year deficit \\ (MWh), 
	ybar=0pt,
	bar width=20pt,
	xtick={1,2,3},
	xticklabels={proposed  method, {MPC,myopic}, {rule-based,myopic}},
	xmin=.5, xmax=3.5
]

\addplot[ao(english), fill=ao(english)!40] table[x=x, y=deficit] {energy_mismatch_method.dat};
\end{axis}
\end{tikzpicture}
  \caption{End-of-year cumulative deficit under three algorithms. Under the proposed framework, every consumer receives at least the amount of promised energy. The proposed MPC + myopic settlement deliver a 0.42 MWh deficit. The deficit with rule-based control + myopic settlement is 1.2 MWh.}
  
  \label{fig:enerfy_deficit_comparison}

\end{figure}
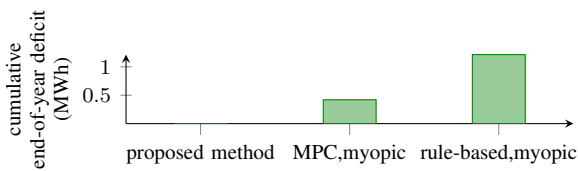
	Fig.~\ref{fig:energy_mismatch}(a) shows that the end-of-year mismatch for each consumer with our proposed MPC and settlement algorithms is between +3\% and +11\%.  In absolute terms, all consumers finish the year with a surplus of 276 $\pm$2 kWh (just over 2\% of total load). As shown in Fig.~\ref{fig:enerfy_deficit_comparison}, under the proposed method every consumer receives at least as much energy as promised.  Note that the main reason all consumers were allocated more energy than promised during the benefit allocation phase, is that PV production during the simulated year was 6\% higher than the expected. As shown in Figs.~\ref{fig:energy_mismatch}(b) and~\ref{fig:energy_mismatch}(c), the alternatives deliver a less equitable energy distribution. In the MPC + myopic settlement case, the energy mismatch per consumer ranges from -4\% to +12\% (the cumulative deficit in this case is 0.42 MWh as shown in Fig.~\ref{fig:enerfy_deficit_comparison}). In the rule-based control + myopic settlement case, the mismatches range from -15\% to +20\% (the cumulative deficit in this case is 1.2 MWh as shown in Fig.~\ref{fig:enerfy_deficit_comparison}). Under the alternatives to our proposed algorithm, some consumers receive less energy than promised (despite having a better-than-average year in terms of PV generation) while others receive over 20\% more than promised. 
	
	The two optimization problems used in this stage are the control problem from the MPC algorithm and the ex-post operation settlement problem. Both are solved at each iteration of MPC. Both problems are convex and easy to solve with commercial solvers. In our case, the computational effort needed to solve these problems is minimal: the control problem solves on an average of $\approx$0.006 seconds, and the settlement problem in $\approx$0.005.

	\begin{figure}
\centering
\begin{tikzpicture}
\begin{axis}[    ylabel= mismatch (\%),
            xlabel= PV capacity factor (p.u.),
                          y label style={at={(axis description cs:-0.1,.5)},rotate=90,anchor=south,align=center, font=\footnotesize},
                x label style={at={(axis description cs:0.5,-.15)},anchor=north, font=\footnotesize},
              width=0.9\linewidth, height=4cm,
              legend style={at={(1,0.3)},
         anchor=east,legend columns=1,draw=none, font=\footnotesize}, 
         xmin = .18, xmax= .20,
         xtick={.18, .19,  .2}, 
         x tick label style={
        /pgf/number format/.cd,
            fixed,
            fixed zerofill,
            precision=2,
        /tikz/.cd
    },
   axis x line*=bottom,
]
\addplot[mark=none, gray, dashed] coordinates {(0.186480491,-10) (0.186480491,10)};
\addplot [only marks, color=blue] table[x= solar_cf, y=cons1] {solar_cf_mismatch.dat};
\addplot [ only marks, color=brown] table[x= solar_cf, y=cons2] {solar_cf_mismatch.dat};
\addplot [ only marks, color=gray] table[x= solar_cf, y=cons3] {solar_cf_mismatch.dat};
\addplot [ only marks, color=green] table[x= solar_cf, y=cons4] {solar_cf_mismatch.dat};
\addplot [ only marks, color=lightgray] table[x= solar_cf, y=cons5] {solar_cf_mismatch.dat};
\addplot [ only marks, color=lime] table[x= solar_cf, y=cons6] {solar_cf_mismatch.dat};
\addplot [ only marks, color=magenta] table[x= solar_cf, y=cons7] {solar_cf_mismatch.dat};
\addplot [ only marks, color=olive] table[x= solar_cf, y=cons8] {solar_cf_mismatch.dat};
\addplot [ only marks, color=orange] table[x= solar_cf, y=cons9] {solar_cf_mismatch.dat};
\addplot [ only marks, color=pink] table[x= solar_cf, y=cons10] {solar_cf_mismatch.dat};
\addplot [ only marks, color=purple] table[x= solar_cf, y=cons11] {solar_cf_mismatch.dat};
\addplot [ only marks, color=red] table[x= solar_cf, y=cons12] {solar_cf_mismatch.dat};
\addplot [ only marks, color=teal] table[x= solar_cf, y=cons13] {solar_cf_mismatch.dat};
\addplot [ only marks, color=violet] table[x= solar_cf, y=cons14] {solar_cf_mismatch.dat};
\addplot [ only marks, color=yellow] table[x= solar_cf, y=cons15] {solar_cf_mismatch.dat};

 \legend{expected PV capacity factor}
\end{axis}
\end{tikzpicture}
\caption{End-of-year mismatch as a function of PV capacity factor for each consumer. The capacity factor corresponds to a realized PV generation scenario.  Each mark represents each of the 15 consumers. The dashed line is the expected PV capacity factors.}

\label{fig:mismatch_CF}
\end{figure}
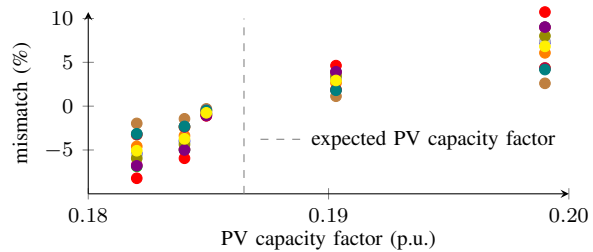

	
	Naturally, the energy allocation mismatch is a function of PV generation. As shown in Fig.~\ref{fig:energy_mismatch}, our algorithm over-delivers PV+S in the simulated year mostly due to a 6\% better-than-average year in terms of PV generation. Fig.~\ref{fig:mismatch_CF} shows the relationship between PV capacity factor (i.e., average generation over installed capacity) and end-of-year mismatch for each consumer for five simulated years. As expected, the allocation surplus (i.e., positive mismatch) tends to increase with capacity factor (i.e., more PV generation). Note that by linear interpolation, at expected capacity factor (denoted by the dashed line in Fig.~\ref{fig:mismatch_CF}) the expected mismatch is close to zero. 
	
	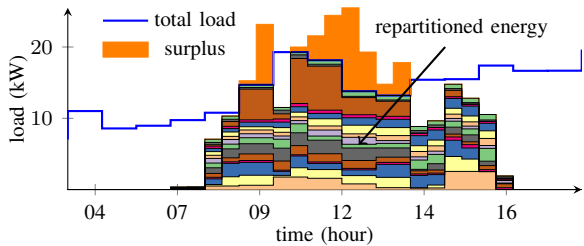
\begin{figure}
\centering
\begin{tikzpicture}

\begin{axis}[
       cycle list/Accent,
  date coordinates in=x,
  ticklabel style = {font=\footnotesize},
  xticklabel=\hour,
  xmin  = 00-01-00 04:00, xmax = 00-01-00 19:00,
  width=.95\linewidth, height=4cm, 
  y label style={at={(axis description cs:-0.06,.5)},rotate=90,anchor=south,align=center, font=\footnotesize},
    x label style={at={(axis description cs:0.5,-0.15)},anchor=north, font=\footnotesize},
	ylabel=load (kW), xlabel= time (hour),
		legend style={at={(0.05,0.85)}, draw=none,fill=none,
		anchor=west,legend columns=1, font=\footnotesize},
]
\path[name path=axis] (axis cs:00-01-00 00:00,0) -- (axis cs:00-01-00 23:30,0);
\addplot[const plot,name path=LOAD, color=blue, thick] table [x=time,y=load, col sep=comma] {actual_key_sample.dat};
\addplot[const plot, name path=SOLAR, color=orange, line width=0pt,forget plot] table [x=time,y=production, col sep=comma] {actual_key_sample.dat};
\addplot[const plot,name path=B, color=black, line width=0pt,forget plot] table [x=time,y =x1 , col sep=comma] {actual_key_sample.dat};
\addplot[const plot,name path=C, color=black, line width=0pt,forget plot] table [x=time,y =x2 , col sep=comma] {actual_key_sample.dat};
\addplot[const plot,name path=D, color=black, line width=0pt,forget plot] table [x=time,y =x3 , col sep=comma] {actual_key_sample.dat};
\addplot[const plot,name path=E, color=black, line width=0pt,forget plot] table [x=time,y =x4 , col sep=comma] {actual_key_sample.dat};
\addplot[const plot,name path=F, color=black, line width=0pt,forget plot] table [x=time,y =x5 , col sep=comma] {actual_key_sample.dat};
\addplot[const plot,name path=G, color=black, line width=0pt,forget plot] table [x=time,y =x6 , col sep=comma] {actual_key_sample.dat};
\addplot[const plot,name path=H, color=black, line width=0pt,forget plot] table [x=time,y =x7 , col sep=comma] {actual_key_sample.dat};
\addplot[const plot,name path=I, color=black, line width=0pt,forget plot] table [x=time,y =x8 , col sep=comma] {actual_key_sample.dat};
\addplot[const plot,name path=J, color=black, line width=0pt,forget plot] table [x=time,y =x9 , col sep=comma] {actual_key_sample.dat};
\addplot[const plot,name path=K, color=black, line width=0pt,forget plot] table [x=time,y =x10 , col sep=comma] {actual_key_sample.dat};
\addplot[const plot,name path=L, color=black, line width=0pt,forget plot] table [x=time,y =x11 , col sep=comma] {actual_key_sample.dat};
\addplot[const plot,name path=M, color=black, line width=0pt,forget plot] table [x=time,y =x12 , col sep=comma] {actual_key_sample.dat};
\addplot[const plot,name path=N, color=black, line width=0pt,forget plot] table [x=time,y =x13 , col sep=comma] {actual_key_sample.dat};
\addplot[const plot,name path=O, color=black, line width=0pt,forget plot] table [x=time,y =x14 , col sep=comma] {actual_key_sample.dat};
\addplot[const plot,name path=P, color=black, line width=0pt,forget plot] table [x=time,y =x15 , col sep=comma] {actual_key_sample.dat};
\addplot[color=orange] fill between[of=P and SOLAR];
\addplot fill between[of=B and axis,forget plot];
\addplot fill between[of=C and B,forget plot];
\addplot fill between[of=D and C,forget plot];
\addplot fill between[of=E and D,forget plot];
\addplot fill between[of=F and E,forget plot];
\addplot fill between[of=G and F,forget plot];
\addplot fill between[of=H and G,forget plot];
\addplot fill between[of=I and H,forget plot];
\addplot fill between[of=J and I,forget plot];
\addplot fill between[of=K and J,forget plot];
\addplot fill between[of=L and K,forget plot];
\addplot fill between[of=M and L,forget plot];
\addplot fill between[of=N and M,forget plot];
\addplot fill between[of=O and N,forget plot];
\addplot fill between[of=P and O,forget plot];
\legend{total load, surplus}

\node[anchor=south] (source) at (axis cs:00-01-00 15:30, 20){\footnotesize repartitioned energy};
       \node (destination) at (axis cs: 00-01-00 12:15,5){};
       \draw[->, thick](source)--(destination);
\end{axis}

\end{tikzpicture}

  \caption{Total load, the key of repartition, and surplus for a sample day. The repartitioned energy is represented by the colored area under the blue curve with each color representing a different resident.}
  \label{fig:key_example}
\end{figure}
	
	Fig.~\ref{fig:key_example} illustrates the allocation of energy in the collective for one sample day. Each of the colored areas under the load represents the energy allocated to a consumer on a 30-minute basis, i.e., the key of repartition. The orange area over the load represents the surplus energy sold by the investor to the grid. Since we present results for the pessimistic case, the surplus goes uncompensated.

	\section{Conclusions and Future Work}
	\label{sec:Conclusions}
	In this paper, we present a framework that integrates long-term planning and operation of a collective that shares a photovoltaic plus storage (PV+S) system. In our framework, the collective is constituted by an investor who provides capital for the PV+S system and a set of consumers who buy PV+S energy at an agreed price from the investor and complement their consumption from the grid. While the presented framework is fairly generic, we focus on the collective self-consumption scheme laid out by French regulation.  
	
	In the long-term planning stage, we first size the PV+S system by solving a mixed-integer linear program whose objective is to maximize the long-term (e.g., 20 years) welfare of all participants of the collective (i.e., the investor and the consumers) and then allocate the expected benefits among the participants. For the benefit allocation stage, we first find the set of PV+S energy prices that would lead to a financially stable collective, i.e., the set of prices that are high enough for the investor to profit but low enough to be attractive for the consumers. Finally, we allocate PV+S among the consumers based on the principle of equity. We accomplish this, by solving an optimization problem whose objective is to allocate PV+S energy evenly among consumers.  
	
	The long-term planning and the operation stages are coupled by the sizing decisions, as usual, but also by the expected PV+S energy allotment to each consumer. Thus, the objective of the operation is both to minimize system operation and to minimize the mismatch between the expected and actual PV+S allotment. The operation stage is divided into two steps. First, we determine the operation strategy on a rolling horizon via a model-predictive control algorithm. Then, after the strategy is implemented and the uncertainty materializes, we adjust the energy repartition to reflect the actual operation of the system.
	
	We demonstrate our framework with a case study of a potential 15 consumer collective in the south of France. We examine two cases: a baseline case that features PV investment costs on the low-end and reasonable prices for surplus sales to the grid and a pessimistic case that features high-end PV investment costs and no remuneration for injections to the grid. We show that in both cases, installing PV+S is Pareto-optimal. Then we show the range of PV+S prices that fosters financially sustainable (prices between 0.10 and 0.13EUR/kWh in the pessimistic case). Then, we show that, even in the pessimistic case and when the investor extracts half of the benefits, the consumer can still expect to save about 3 to 5\%. Naturally, the savings would be higher with more favorable consumer conditions (e.g., lower PV installation costs or higher share of benefits assigned to consumers). Finally, our proposed control and settlement algorithms outperform two alternatives in terms of the end-of-year mismatch between energy promised and energy delivered.  
	
	Future work can take diverse directions. For example, we did not consider consumers equipped with demand response (DR). It would be interesting to investigate the effects of DR in different indices such as system cost and rate of self-consumption and problems such as DR control and compensation. Furthermore, our model mainly deals with the energy management problem. One could enhance the model by including the local distribution grid and considering other voltage stability or power quality problems.

	\ifCLASSOPTIONcaptionsoff
	\newpage
	\fi

	
	
	
	
	\bibliographystyle{IEEEtran}
	\bibliography{Bibliography}

\begin{thebibliography}{10}
\providecommand{\url}[1]{#1}
\csname url@rmstyle\endcsname
\providecommand{\newblock}{\relax}
\providecommand{\bibinfo}[2]{#2}
\providecommand\BIBentrySTDinterwordspacing{\spaceskip=0pt\relax}
\providecommand\BIBentryALTinterwordstretchfactor{4}
\providecommand\BIBentryALTinterwordspacing{\spaceskip=\fontdimen2\font plus
\BIBentryALTinterwordstretchfactor\fontdimen3\font minus
  \fontdimen4\font\relax}
\providecommand\BIBforeignlanguage[2]{{%
\expandafter\ifx\csname l@#1\endcsname\relax
\typeout{** WARNING: IEEEtran.bst: No hyphenation pattern has been}%
\typeout{** loaded for the language `#1'. Using the pattern for}%
\typeout{** the default language instead.}%
\else
\language=\csname l@#1\endcsname
\fi
#2}}

\bibitem{eu2015best}
{European Union Commission}, ``Best practices on renewable energy
  self-consumption,'' \emph{Commission staff working document}, 2015.

\bibitem{Frieden_2019_collective}
D.~Frieden, A.~Tuerk, J.~Roberts, S.~d'Herbemont, and A.~Gubina, ``{Collective
  self-consumption and energy communities: Overview of emerging regulatory
  approaches in Europe},'' June 2019, working paper, Compile project.

\bibitem{code_energie}
\BIBentryALTinterwordspacing
``{Code de l'\'energie}.'' [Online]. Available:
  \url{https://www.legifrance.gouv.fr/}
\BIBentrySTDinterwordspacing

\bibitem{deliberation_CRE}
{Commission de R\'egulation de l'\'Energie}, ``{D\'ELIBERATION No. 2018-115},''
  2018.

\bibitem{abada2017viability}
I.~Abada, A.~Ehrenmann, and X.~Lambin, ``On the viability of energy
  communities,'' \emph{Energy Journal}, 2019.

\bibitem{Dunlop_2016_EU}
S.~Dunlop and A.~Roesch, ``{EU-Wide Solar PV Business Models},'' November 2016,
  {SolarPower Europe}.

\bibitem{goedkoop2016partnership}
F.~Goedkoop and P.~Devine-Wright, ``{Partnership or placation? The role of
  trust and justice in the shared ownership of renewable energy projects},''
  \emph{Energy Research \& Soc. Science}, vol.~17, pp. 135--146, 2016.

\bibitem{shared_schroeder}
\BIBentryALTinterwordspacing
E.~M. Schroeder, ``{Shared Solar Programs: Opportunities and Challenges},''
  presentation. [Online]. Available:
  \url{https://nrel.gov/state-local-tribal/assets/pdfs/stat_webinar_070913_presentation.pdf}
\BIBentrySTDinterwordspacing

\bibitem{Council_2019_regulatory}
{Council of European Energy Regulators}, ``{Regulatory Aspects of
  Self-Consumption and Energy Communities},'' 2019, c18-CRM9\_DS7-05-03.

\bibitem{EU_RE_directive}
{The European Union Parliament and The Council of The European Union},
  ``Directive (eu) 2018/2001 on the promotion of the use of energy from
  renewable sources,'' 2018,
  \newline\url{https://eur-lex.europa.eu/legal-content/EN/TXT/?uri=uriserv:OJ.L_.2018.328.01.0082.01.ENG&toc=OJ:L:2018:328:TOC}.

\bibitem{EU_IME_directive}
------, ``Directive (eu) 2018/2001 on common rules for the internal market for
  electricity and amending directive 2012/27/eu,'' 2019,
  \newline\url{https://eur-lex.europa.eu/legal-content/EN/TXT/?uri=CELEX\%3A32019L0944}.

\bibitem{Frieden_2018_Overview}
D.~{Frieden}, J.~{Roberts}, and A.~F. {Gubina}, ``Overview of emerging
  regulatory frameworks on collective self-consumption and energy communities
  in europe,'' in \emph{2019 16th International Conference on the European
  Energy Market (EEM)}, Sep. 2019, pp. 1--6.

\bibitem{alaton2020energy}
C.~Alaton, J.~Contreras-Oca{\~n}a, P.~de~Radigu{\`e}s, T.~D{\"o}ring, and
  F.~Tounquet, ``Energy communities: From european law to numerical modeling,''
  \emph{preprint arXiv:2008.03044}, 2020.

\bibitem{Moret_2019_Energy}
F.~{Moret} and P.~{Pinson}, ``Energy collectives: A community and fairness
  based approach to future electricity markets,'' \emph{{IEEE} Trans. Power
  Syst.}, vol.~34, no.~5, pp. 3994--4004, Sep. 2019.

\bibitem{Vanadzina_2019_Business}
E.~{Vanadzina}, G.~{Mendes}, S.~{Honkapuro}, A.~{Pinomaa}, and H.~{Melkas},
  ``Business models for community microgrids,'' in \emph{2019 16th
  International Conference on the European Energy Market (EEM)}, Sep. 2019.

\bibitem{Fleischhacker_2019_Sharing}
A.~{Fleischhacker}, H.~{Auer}, G.~{Lettner}, and A.~{Botterud}, ``Sharing solar
  pv and energy storage in apartment buildings: Resource allocation and
  pricing,'' \emph{{IEEE} Trans. Smart Grid}, vol.~10, no.~4, pp. 3963--3973,
  July 2019.

\bibitem{Nunna_2017_Multiagent}
H.~S. V. S.~K. {Nunna} and D.~{Srinivasan}, ``Multiagent-based transactive
  energy framework for distribution systems with smart microgrids,''
  \emph{{IEEE} Trans. Indust. Inf.}, vol.~13, no.~5, pp. 2241--2250, Oct 2017.

\bibitem{Park_2017_Event}
S.~{Park}, J.~{Lee}, G.~{Hwang}, and J.~K. {Choi}, ``Event-driven energy
  trading system in microgrids: Aperiodic market model analysis with a game
  theoretic approach,'' \emph{IEEE Access}, vol.~5, pp. 26\,291--26\,302, 2017.

\bibitem{Shamshi_2016_Shamsi}
P.~{Shamsi}, H.~{Xie}, A.~{Longe}, and J.~{Joo}, ``Economic dispatch for an
  agent-based community microgrid,'' \emph{{IEEE} Trans. Smart Grid}, vol.~7,
  no.~5, pp. 2317--2324, Sep. 2016.

\bibitem{Liu_2019_ANovel}
W.~{Liu}, J.~{Zhan}, and C.~Y. {Chung}, ``A novel transactive energy control
  mechanism for collaborative networked microgrids,'' \emph{{IEEE} Trans. Power
  Syst.}, vol.~34, no.~3, pp. 2048--2060, May 2019.

\bibitem{Feng_2019_Caolitional}
C.~{Feng}, F.~{Wen}, S.~{You}, Z.~{Li}, F.~{Shahnia}, and M.~{Shahidehpour},
  ``Coalitional game based transactive energy management in local energy
  communities,'' \emph{{IEEE} Trans. Power Syst.}, pp. 1--1, 2019.

\bibitem{Etinski_2013_Fair}
M.~{Etinski} and A.~{Sch\"ulke}, ``Fair power allocation in multi-user systems
  with controllable loads,'' in \emph{2013 IEEE SmartGridComm}, Oct 2013, pp.
  37--42.

\bibitem{Chen_2019_PV}
X.~{Chen}, J.~{Zhao}, and M.~{He}, ``{PV} power generation credit sharing
  towards sustainable community solar,'' in \emph{2019 IEEE TPEC}, Feb 2019,
  pp. 1--6.

\bibitem{Khare_2014_Optimal}
A.~{Khare} and S.~{Rangnekar}, ``Optimal sizing of a grid integrated solar
  photovoltaic system,'' \emph{IET Renew. Power Gen.}, vol.~8, no.~1, pp.
  67--75, January 2014.

\bibitem{Kolhe_2009_Techno}
M.~{Kolhe}, ``Techno-economic optimum sizing of a stand-alone solar
  photovoltaic system,'' \emph{{IEEE} Trans. Energy Convers.}, vol.~24, no.~2,
  pp. 511--519, June 2009.

\bibitem{Contreras-Ocana_2019_Non-Wire}
J.~E. {Contreras-Oca\~na}, Y.~{Chen}, U.~{Siddiqi}, and B.~{Zhang}, ``Non-wire
  alternatives: An additional value stream for distributed energy resources,''
  \emph{{IEEE} Trans. Sust. Energy}, pp. 1--1, 2019.

\bibitem{Shahidehpour_2019_Optimal}
M.~{Shahidehpour}, C.~{Li}, H.~{Yang}, B.~{Zhou}, Y.~{Cao}, L.~{Zeng}, and
  Z.~{Xu}, ``Optimal planning of islanded integrated energy system with
  solar-biogas energy supply,'' pp. 1--1, 2019.

\bibitem{Xu_2013_AnImporved}
L.~{Xu}, X.~{Ruan}, C.~{Mao}, B.~{Zhang}, and Y.~{Luo}, ``An improved optimal
  sizing method for wind-solar-battery hybrid power system,'' \emph{{IEEE}
  Trans. Sust. Energy}, vol.~4, no.~3, pp. 774--785, July 2013.

\bibitem{Akram_2018_AnImproved}
U.~{Akram}, M.~{Khalid}, and S.~{Shafiq}, ``An improved optimal sizing
  methodology for future autonomous residential smart power systems,''
  \emph{IEEE Access}, vol.~6, pp. 5986--6000, 2018.

\bibitem{Saez-de-Ibarra_2019_Co-Optimization}
A.~{Saez-de-Ibarra}, A.~{Milo}, H.~{Gazta\~naga}, V.~{Debusschere}, and
  S.~{Bacha}, ``Co-optimization of storage system sizing and control strategy
  for intelligent photovoltaic power plants market integration,'' \emph{{IEEE}
  Trans. Sust. Energy}, vol.~7, no.~4, pp. 1749--1761, Oct 2016.

\bibitem{levy_1994_capital}
H.~Levy and M.~Sarnat, \emph{Capital investment and financial decisions}.\hskip
  1em plus 0.5em minus 0.4em\relax Pearson Education, 1994.

\bibitem{Hannoset_2019_Energy}
A.~Hannoset, L.~Peeters, and A.~Tuerk, ``{Energy Communities in the EU Task
  Force Energy Communities},'' 2019, {BRIDGE H2020 Initiative}.

\bibitem{Contreras_2019_Participation}
J.~E. {Contreras-Oca\~na}, M.~A. {Ortega-Vazquez}, and B.~{Zhang},
  ``Participation of an energy storage aggregator in electricity markets,''
  \emph{{IEEE} Trans. Smart Grid}, vol.~10, no.~2, pp. 1171--1183, March 2019.

\bibitem{Sarker_2017_Optimal}
M.~R. Sarker, M.~D. Murbach, D.~T. Schwartz, and M.~A. Ortega-Vazquez,
  ``Optimal operation of a battery energy storage system: Trade-off between
  grid economics and storage health,'' \emph{Electric Power Systems Research},
  vol. 152, pp. 342 -- 349, 2017.

\bibitem{nash1950bargaining}
J.~F. Nash~Jr, ``The bargaining problem,'' \emph{Econometrica: Journal of the
  Econometric Society}, pp. 155--162, 1950.

\bibitem{binmore1986nash}
K.~Binmore, A.~Rubinstein, and A.~Wolinsky, ``The nash bargaining solution in
  economic modelling,'' \emph{The RAND Journal of Economics}, pp. 176--188,
  1986.

\bibitem{Singh_2019_Modeling}
A.~Singh, ``Modeling of a collective self-consumption system in a residential
  district,'' 2019, {Final year internship report. Grenoble INP, Arkolia
  Energies}.

\bibitem{jensen_RE-Europe_data}
\BIBentryALTinterwordspacing
T.~V. Jensen, H.~de~Sevin, M.~Greiner, and P.~Pinson, ``{The RE-Europe data
  set},'' 2015. [Online]. Available: \url{https://doi.org/10.5281/zenodo.35177}
\BIBentrySTDinterwordspacing

\bibitem{pecan_street}
``{Pecan Street Dataport},'' \url{https://dataport.pecanstreet.org/}, accessed:
  2019-09-19.

\bibitem{photovoltaique.com}
\BIBentryALTinterwordspacing
``Le centre de ressources photovolta\"ique.'' [Online]. Available:
  \url{https://www.photovoltaique.info/fr/}
\BIBentrySTDinterwordspacing

\bibitem{edf_france_2018}
\BIBentryALTinterwordspacing
``{Quelles aides \`a l'autoconsommation pour les collectivit\'es? - EDF
  Collectivit\'es},'' Jun 2018. [Online]. Available:
  \url{https://www.edf.fr/en/node/295857}
\BIBentrySTDinterwordspacing

\bibitem{fu2018us}
R.~Fu, R.~M. Margolis, and D.~J. Feldman, ``{US Solar Photovoltaic System Cost
  Benchmark: Q1 2018},'' National Renewable Energy Lab.(NREL), Golden, CO
  (United States), Tech. Rep., 2018.

\bibitem{Lavergne_2019_Stockage}
R.~Lavergne, I.~Pavel, and I.~Faucheux, ``{Stockage Stationnaire
  d'\'Electricit\'e},'' Mar 2019, {Conseil G\'en\'eral de l'\'Economie}.

\bibitem{weniger_2014_sizing}
J.~Weniger, T.~Tjaden, and V.~Quaschning, ``{Sizing of Residential PV Battery
  Systems},'' \emph{Energy Procedia}, vol.~46, pp. 78 -- 87, 2014, 8th Intl.
  Renewable Energy Storage Conf. and Exhibition.

\end{thebibliography}

	\vfill
	

\end{document}